\documentclass[12pt]{article}

\usepackage{amsthm,amsmath,amsfonts}
\usepackage{latexsym}
\usepackage[T1]{fontenc}
\usepackage{a4wide}
\usepackage{times}
\usepackage{array}
\usepackage{hhline}
\usepackage{graphicx}
\usepackage{color}

\numberwithin{equation}{section}

\newtheorem{tw}{Theorem}
\newtheorem{lem}{Lemma}

\begin{document}
\begin{center}
{\bf \large Compressible perturbation of Poiseuille type flow}
\vskip5mm
Piotr B. Mucha ~~ and ~~ Tomasz Piasecki

\vskip3mm

{Institute of Applied Mathematics and Mechanics, University of Warsaw}

{ul. Banacha 2, 02-097 Warszawa, Poland }

\vskip3mm

{E-mail: {\tt p.mucha@mimuw.edu.pl}}

\vskip3mm

{E-mail: {\tt tpiasecki@mimuw.edu.pl}}

\end{center}
\vskip5mm

{\bf Abstract.} The paper examines the issue of stability of Poiseuille type flows in regime of 
compressible Navier-Stokes equations in a three dimensional finite pipe-like domain. 
We prove the existence of stationary solutions with inhomogeneous Navier slip boundary conditions admitting
nontrivial inflow condition in the vicinity of constructed generic flows. Our techniques are based on an application of a modification 
of the Lagrangian coordinates. Thanks to such approach we are able to overcome difficulties 
coming from hyperbolicity of the continuity equation, constructing a maximal regularity estimate for a linearized system 
and applying the Banach fixed point theorem.

\smallskip 



\section{Introduction}

The mathematical description of compressible flows is important from the point of view of applications, 
domains such as aerodynamics and geophysics are the most natural to be mentioned here. 
On the other hand, complexity of the equations describing the flow delivers 
very interesting mathematical challenges.  
In spite of active research in the field, 
we are still far from the complete mathematical understanding of compressible flows. 
The only general existence results are available for weak solutions with homogeneous boundary conditions \cite{Fe}, \cite{PL}.
As far as regular solutions are concerned, we have so far only partial results
assuming either some smallness of the data, or its special structure. The problems have been investigated mainly with
homogeneous boundary conditions (\cite{BdV}, \cite{NoPa}).
For the overview of the state of art in the theory one can consult the monograph \cite{NoS}.

From the point of view of the aforementioned applications 
it seems very important to investigate the problems with large 
velocity vectors, which lead in a natural way to inhomogeneous boundary conditions.
Due to the hyperbolic character of the continuity equation the density must be then 
prescribed on the inflow part of the boundary.
Existence issues for such inflow problems are investigated in \cite{Kw1}, \cite{Kw2}, \cite{TP1}, \cite{TP2},  \cite{PRS1},  \cite{PRS2}
and \cite{VZ}. 
%
%
The mentioned group of problems can be regarded as questions of stability of particular constant flows. 
 
In the present article we would like to examine the issue of stability of Poiseuille type 
flow in pipe-like domain in compressible regime. The unperturbed flow is a solution to the 
compressible Navier-Stokes system with homogeneous slip boundary conditions for given constant 
friction, constant density and constant external force (gravitation-like term). 
The Poiseuille flow is a special 
symmetric solution to the incompressible Navier-Stokes equations in cylindrical domains. Here 
it is viewed as a solution to the compressible Navier-Stokes system with constant density
and constant external force, parallel to axis of the cylinder, given by the pressure.
Hence the pressure, unknown in the incompressible model, is recognized as a given force. 
Such change of `observer' looks acceptable from the mechanical point of view. 
Thanks to that interpretation we obtain a natural physically reasonable flow in compressible
regime. The mathematical objective of this article is to establish stability 
of such flow under some structural assumptions limiting the magnitude of admissible perturbations.

\begin{figure}[htb]
\begin{center}
\includegraphics[width = 0.5\textwidth]{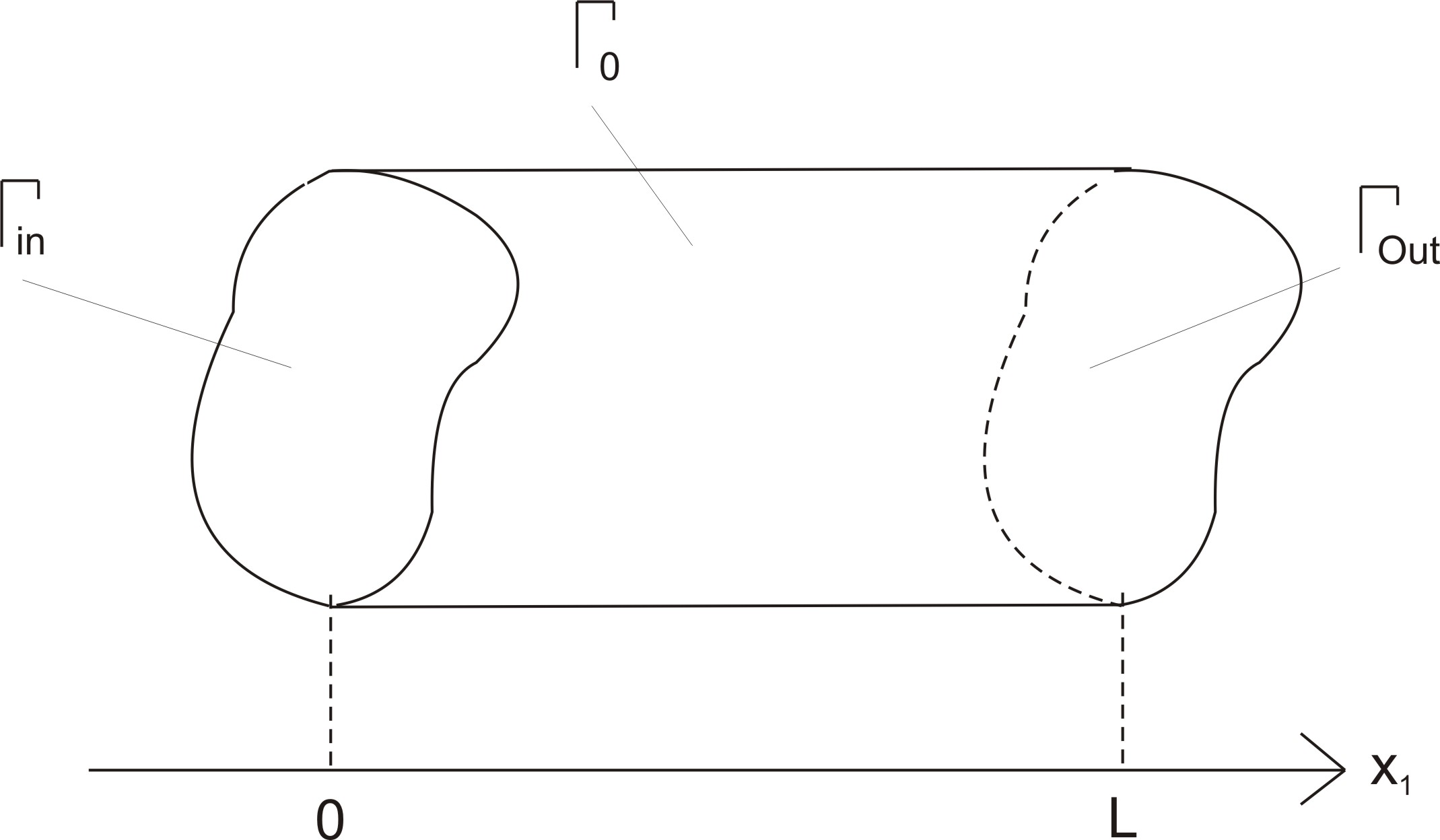}
\caption{The domain}      \label{rys1}
\end{center}
\end{figure}

Let us define the system.
We consider steady flow of a viscous, barotropic fluid in a bounded, cylindrical
domain in $\mathbb{R}^{3}$, described by the Navier - Stokes system
supplied with inhomogeneous Navier slip boundary conditions. The complete system reads
\begin{equation}  \label{main_system}
\begin{array}{lcr}
\rho v \cdot \nabla v -\mu \Delta v - (\mu+\nu) \nabla {\rm div}\; v
+\nabla \pi(\rho) = \rho \, F & \mbox{in} & \Omega,\\
{\rm div}\;(\rho v)=0 & \mbox{in} & \Omega,\\
n\cdot {\bf T}(v,\pi)\cdot \tau_k +f v\cdot \tau_k=b_k, \quad k=1,2 &
\mbox{on} &\Gamma,\\
n \cdot v=d & \mbox{on} & \Gamma,\\
\rho=\rho_{in} & \mbox{on} & \Gamma_{in},
\end{array}
\end{equation}
where $\Omega = [0,L] \times \Omega_0$ with $\Omega_0 \subset \mathbb{R}^2$ of class $C^2$, $\Gamma$ denotes the boundary of $\Omega$ (see Fig.1),
$v$ is the velocity field of the fluid, 
$\rho$ its density, $\mu$ and $\nu$ are viscosity constants satisfying $\mu>0$ and
$(\nu + 2 \mu) >0$, $f \geq 0$ is the friction coefficient which may be different on different components 
of the boundary $\Gamma$, $\pi=\pi(\rho)$ is the pressure given as a function, at least $C^{1}$, of the density and $F$ is an external force.
${\bf T}$ denotes the Cauchy stress tensor of the form 
${\bf T}(v,\pi)=2 \mu {\bf D}v +\nu {\rm div} v {\bf Id} - \pi {\bf Id}$ 
where ${\bf D} = \frac{1}{2}(\nabla v+\nabla v^T)$ is the symmetric gradient. 
Next, $n$ and $\tau_k$ are outer normal and tangent vectors to $\partial \Omega$. 
Boundary data $\rho_{in},b,d$ will be discussed later. The boundary $\Gamma$ is naturally split into three parts:
\begin{equation}
\begin{array}{c}
\Gamma_0=\{x\in \partial\Omega: v(x)\cdot n(x) =0\}, \\
\Gamma_{in}=\{x\in \partial\Omega: v(x)\cdot n(x) < 0\}, \\
\Gamma_{out}=\{x\in \partial\Omega: v(x)\cdot n(x) >0\}.
\end{array}
\end{equation}
Thanks to the chosen geometry of the domain, the above decomposition is easily illustrated by the figure \ref{rys1}.

We shall say few words about the physical interpretation of the system (\ref{main_system}),
in particular about the choice of boundary conditions (\ref{main_system})$_{3,4}$.
We would like to model a flow through a pipe. 
We assume that the fluid obeys Navier slip conditions on the walls of the pipe 
($\Gamma_0$ component of the boundary), hence natural conditions on $\Gamma_0$ are $d \equiv 0$
and $b_k \equiv 0$. However, the mathematical requirements impose a need to prescribe the boundary
conditions on $\Gamma_{in}$ and $\Gamma_{out}$. From the physical viewpoint these parts are artificial,
this is the area where the parameters of the velocity and density are measured. This gives us
a freedom of choice of the type of boundary conditions on the inflow and outflow part, which can be fit
to the mathematical approach, hence
we choose inhomogeneous slip condition.  
Note that as the friction coefficient goes to infinity, then
the relations $(\ref{main_system})_{3,4}$, at least formally, become the standard Dirichlet conditions 
describing the whole velocity vector at the boundary.
Since the velocity does not vanish on the boundary, the hyperbolicity of the continuity
equation impose a need to prescribe the density on the inflow part, which lead to the condition
(\ref{main_system})$_5$. The velocity field determines the characteristics of the continuity equation and in particular
the total mass $\int_{\Omega}\rho \,dx $ is determined implicitly by $(\ref{main_system})_5$.

Our goal here is to analyze
a perturbation of the Poiseuille type flow 
\begin{equation} \label{Pois_gen}
\bar V = [V^P(x_2,x_3),0,0],
\end{equation}  
where $x_1$ points the direction of the axis of the cylinder.
It is one of the classical examples of laminar flows satisfying the incompressible
Navier-Stokes equations in cylindrical domains. 
In the classical literature the flow is considered with homogeneous Dirichlet condition on the boundary.
$V^P$ is then found as a solution to corresponding elliptic problem with Dirichlet boundary condition 
on each $x_1$ - cut of $\Omega$. Some explicit formulas on $V^P$ in certain domains
are well-known (\cite{Ga},\cite{LadSol}).

In the case of slip boundary conditions that are subject of our analysis in this paper 
the flow (\ref{Pois_gen})
can be also found on each cut of the cylinder as a solution to elliptic problem with corresponding
boundary condition (see Lemma 1 below).
In certain domains it can also be expressed with explicit formulas 
(see \cite{PM1} and the example below). 
Since we are interested in a general cylindrical domain, we will not have such formula but we  show that the solution of the form (\ref{Pois_gen}) exists         
provided that $\Omega_0$ is sufficiently regular and, under the slip boundary conditions (\ref{main_system})$_{3,4}$,
$V^P$ does not vanish on the boundary.

\begin{lem} \label{lem_vp} 
Let $\Omega_\infty = \mathbb R \times \Omega_0$, where $\Omega_0 \subset \mathbb{R}^2$ is smooth enough,
$\mu>0$ and $f \geq 0$.
Then there exists a solution $(\bar V , \bar \Pi)$, such that $\bar V=(V^P(x_2,x_3),0,0)$,
to the incompressible Navier-Stokes system with
slip boundary conditions:
\begin{equation}
\begin{array}{lcr} \label{eqn_pois_omega}
\bar V \cdot \nabla \bar V - \mu \Delta \bar V + \nabla \bar \Pi = 0 & \mbox{in} & \Omega_\infty, \\
{\rm div} \, \bar V = 0 & \mbox{in} & \Omega_\infty,\\
n\cdot {\bf T}(\bar v, \bar \Pi)\cdot \tau_k +f \bar V\cdot \tau_k = 0, \quad k=1,2 &
\mbox{on} & \mathbb R \times \partial \Omega_0,\\
n \cdot \bar V = 0 & \mbox{on} & \mathbb R \times \partial \Omega_0.\\
\end{array}
\end{equation}
Moreover, there exist $ \theta=\theta(f,\mu)$ and a continuous function $\bar \omega(f)$ such that
\begin{equation} \label{VP_positive} 
\bar V^{(1)}=V^P \geq \theta > 0 \quad \textrm{in} \quad \bar \Omega_\infty.
\end{equation}
and
\begin{equation} \label{est_VP}
\|\nabla V^p\|_{L_\infty} \leq \frac{\bar \omega(f)}{\mu}.
\end{equation}
In addition if $(\bar V,\bar \Pi)$ is the Poiseuille solution, then $(\lambda \bar V, \lambda \bar \Pi)$, too.
\end{lem}

As we already said, the pressure in the Poiseuille flow can be regarded 
as an external force parallel to the axis
of the cylinder. 
Like in the classical Poiseuille flow, we assume the pressure to be a linear function of
$x_1$. 
Hence it is natural to assume the form $\bar \pi = \omega(f) x_1$,
where $\omega(f)$ is a negative constant (the sign describes the directon of the flow).
Then we see that $V^P$, the first coordinate of $\bar V=(V^P,0,0)$, is a solution to the elliptic problem 
\begin{equation} \label{eqn_pois}
\begin{array}{c}
\mu \Delta V^P = \bar \Pi_{x_1} = \omega(f) < 0 \mbox{ ~~ in~} \Omega_0, \\
\mu \frac{\partial V^P}{\partial n} + f \, V^P = 0 \mbox{ ~~ on ~} \partial \Omega_0.
\end{array}
\end{equation}
Now to prove Lemma \ref{lem_vp} it is enough to apply the maximum principle to the system (\ref{eqn_pois}).
We show the proof in the Appendix, at this stage we should
have a closer look at the dependence $\omega(f)$. 
This dependence can be justified considering the compatibility condition 
for the system (\ref{eqn_pois}),
that reads $\int_{\Omega_0} \omega \,dx = f \int_{\partial \Omega_0} V^P \,d\sigma $.
As we will see, the dependence $\omega(f)$ determines the assumptions
we will have to make on the viscosity.
  
Note that for $\omega=0$ the only solution is $V^P=0$, provided $f \neq 0$. 
For $f=0$ (perfect slip) we should obtain a constant flow, hence we put $\omega=0$ for $f=0$.
The linear structure of (\ref{eqn_pois}) makes $\omega(f)$ continuous. 
Thus, we conclude $\omega(f) \to 0$ for $f \to 0$, and hence by (\ref{est_VP})  
\begin{equation} \label{nabla_pois_small} 
\|\nabla \bar V\|_{L_\infty} \to 0 \; \textrm{~~ for~~ } \; f \to 0.
\end{equation}
This is an important conclusion since in order to show the energy estimate 
we have to control $\nabla \bar V$ with the viscosity, and
so we  allow the viscosity to be low provided that the friction on $\Gamma_0$ is small, what
is a realistic assumption (see also the remarks after the formulation of Theorem 1). 
Let us illustrate the dependence $\omega=\omega(f)$ with the following example.

\smallskip

{\bf Example.} Take $\Omega_0 = B(0,1)\subset \mathbb R^2$ and $\mu =1$. Then, 
due to axial symmetry of the domain,
it is natural to look for 
$V^P = v^f(r)$ where $r = \sqrt{x_2^2+x_3^2}$. 
The boundary condition (\ref{eqn_pois})$_2$ then reads $v^f_r + f v^f|_{r=1} = 0$.
We require that 
$\Delta v = v_{rr} + \frac{1}{r}v_r$ depend only on $f$. Moreover, we expect to
obtain a constant flow for $f=0$ (perfect slip) and classical Poiseuille profile for $f = \infty$.
The above considerations lead to the family of solutions 
$$
v^f(r) = TF \big[ \frac{f+2}{f+k_f} - \frac{f}{f+k_f} r^2 \big],
$$
where $k_f = \frac{(\pi-2)f + 4 \pi}{2}$ and $TF$ is the flux of the flow $v^f$ through $\Omega_0$.

For a perfect slip case $f=0$ we obtain a constant flow $v^0 = \frac{2 \, TF}{k_f}$ 
and for a no-slip case $f=\infty$ we get a classical Poiseuille profile $v^\infty = TF(1-r^2)$.
On the boundary we have $\theta(f) = v^f(1) = \frac{2}{f+k_f}$ what is a strictly positive constant.
Finally,
$$
\omega(f) = \Delta v^f = - \frac{4f}{f+k_f} \mbox{~~~~and~~~~} \nabla v^f \sim \frac{f}{f+k_f}.  
$$     
In particular  $\omega(f)<0$ for $f>0$ and $\omega(0) = 0$. 
 
\smallskip

Before we formulate our main result, we need one observation concerning the boundary conditions.
Note that, since $V^P$ is found on every $x_1$ - cut of $\Omega$, we can impose
the boundary conditions (\ref{eqn_pois_omega})$_{3,4}$ only on $\Gamma_0$. 
On the other hand, in order to define small
perturbations as a solution to (\ref{main_system}) we have to 
measure the distance (in appropriate norms) between the solution to (\ref{main_system})
and the Poiseuille flow $\bar V$. Hence we need to consider the traces of the quantities from the
boundary conditions of (\ref{main_system}) with the function $\bar V$ instead of $v$.    
Since our analysis acts on a finite cylinder we define these traces at the bottoms $\Gamma_{in} \cup \Gamma_{out}$:
$$
\bar b_k|_{\Gamma_{in} \cup \Gamma_{out}} = \textrm{tr}_{\Gamma_{in} \cup \Gamma_{out}} 
[n\cdot {\bf T}(\bar V, \bar \Pi)\cdot \tau_k +f \bar V\cdot \tau_k],
\qquad 
\bar d|_{\Gamma_{in}} = - \textrm{tr}_{\Gamma_{in}} V^P,
\qquad
\bar d|_{\Gamma_{out}} = \textrm{tr}_{\Gamma_{out}} V^P,
$$
where $\textrm{tr}$ denotes the trace operator.
By (\ref{eqn_pois_omega}),
$$
n\cdot {\bf T}(\bar V, \bar \Pi)\cdot \tau_k +f \bar V \cdot \tau_k = 0 \mbox{ ~~ and ~~ } n \cdot \bar V = 0 \mbox{ ~~ at } \Gamma_0.
$$
Hence the construction of Poiseuille flow determines $(k=1,2)$:
\begin{equation}\label{d-2}
\bar b_k = \bar d = 0 \mbox{ ~~ ~~ on ~} \Gamma_0.
\end{equation}
The Poiseuille flow $(\bar V, \bar \Pi)$ has constant density $\bar \rho$, we set $\bar \rho=1$.
Then in the chosen setting $(\bar V, \bar \Pi)$ fulfills the following system
\begin{equation}\label{d-1}
\begin{array}{lcr}
\bar \rho \bar V \cdot \nabla \bar V - \mu \Delta \bar V - (\mu+\nu)\nabla {\rm div} \, \bar V + \nabla \bar \Pi(\bar \rho)= - \bar \rho  \omega(f) \hat{e}_1 & \mbox{in} & \Omega,\\
{\rm div} \,(\bar\rho \bar V)=0  & \mbox{in} & \Omega,\\
n\cdot {\bf T}(\bar V, \bar \Pi)\cdot \tau_k +f \bar V\cdot \tau_k = \bar b_k, \quad k=1,2 &
\mbox{on} & \Gamma,\\
n \cdot \bar V = \bar d & \mbox{on} & \Gamma.\\
\end{array}
\end{equation}
We keep in mind that $\nabla \bar \Pi =\omega(f) \hat{e}_1$ and (\ref{d-2}).

We are now in a position to formulate our main result. 
To this end it is convenient to define the quantity which measures the distance of the data from the 
Poiseuille flow: 
\begin{equation} \label{d0}
D_0 = \|F+\omega(f) \hat{e}_1\|_{L_p(\Omega)} + \|d - \bar d\|_{W^{2-1/p}_p(\Gamma)} + \|b_k-\bar b_k\|_{W^{1-1/p}_p(\Gamma)}
+ \|\rho_{in}-1\|_{W^1_p(\Gamma_{in})}.
\end{equation}
The main result of the paper reads

\smallskip 

\begin{tw}
Assume that the boundary data of (\ref{main_system}) is close to the Poiseuille flow $(\bar V, \bar \Pi)$,
more precisely, let $D_0$ defined above be small enough. Assume that the friction $f$ 
is large enough on $\Gamma_{in}$ and $p>3$.
Assume further that there exists a constant $\kappa>0$ such that 
\begin{equation} \label{mu_large}
\mu > \kappa \, \textrm{min} \{ f|_{\Gamma_0},1 \}.
\end{equation}
Then there exists a solution 
$(v,\rho) \in W^2_p(\Omega) \times W^1_p(\Omega)$ to the system (\ref{main_system}) such that 
\begin{equation} \label{est_main}
\|v - \bar V\|_{W^2_p(\Omega)} + \|\rho - \bar \rho\|_{W^1_p(\Omega)} \leq C(D_0).
\end{equation} 
This solution is unique in the class of small perturbations of $(\bar V, \bar \rho)$.
\end{tw}

Let us make some remarks concerning our main result. The condition on the viscosity (\ref{mu_large})
seems to be a serious constraint, but as we will see from the proofs we just need it
to control the gradient of the Poiseuille flow, what yields this assumption natural
(see also the remark in the proof of Lemma \ref{lem_ene1}). 
We recall that $\nabla \bar V$ depends on the friction $f$ on $\Gamma_0$,
and in particular (\ref{nabla_pois_small}) holds. It follows
that for small values of friction on $\Gamma_0$ it is enough to assume that the viscosity is large enough,
but only compared to the friction. This assumption is reflected in the condition (\ref{mu_large}). 
Theorem 1 admits the case of perfect slip $f|_{\Gamma_0} = 0$, and in such case 
(\ref{mu_large}) reduces to $\mu > 0$, so no lower bound on the viscosity is required. 
In this case there is no bound on the size of $\bar V$. 
However in this case $\bar V $ would be a constant flow. 
We shall recall that the friction at $\Gamma_{in}$ is chosen independently to $f$ at $\Gamma_0$. 
From the point of view of modelling, the data at $\Gamma_{in}$ is given, hence it is important 
to focus the attention at $\Gamma_0$. The assumption $p>3$ is required for the imbedding ${W^1_p} \subset {L_\infty}$ \cite{Ad},
it is required to control the pointwise boundedness of  $\nabla v$ and the density.

Let us explain the main idea of the proof. We will follow an idea of Lagrangian type coordinates \cite{PC}, \cite{Mu}, \cite{MZ}, \cite{Sol2}. The continuity equation
is of hyperbolic type and contains a term $u \cdot \nabla w$ 
(where $u$ and $w$ are perturbations to the velocity and density introduced in the next section),
which makes serious troubles for the issues of existence in case of
inhomogeneous boundary conditions, 
see \cite{Kw1}, \cite{PRS1}, \cite{NoS}, \cite{TP1}, \cite{TP2}. Here we overcome this obstacle by changing the system of coordinates 
in such a way that this term disappears (\ref{change_id}). 
We obtain a more complex system but with structure suitable for an application of the 
Banach fixed point theorem. On the other hand 
our solutions are regular enough, thus we are able to go back to the original system keeping the well posedness of the original model.
Our approach works since we are equipped with the maximal regularity estimate for a linearization 
of the  equations in the Lagrangian coordinates -- Theorem \ref{tw_est_main_lin}. 
This tool gives a complete control of the regularity of solutions.

The rest of the paper is organized as follows. 
In Section 2 we introduce the perturbations as unknown variables obtaining the system (\ref{system}). 
Next we introduce the Lagrangian-type coordinates that lead to the system (\ref{system_z})
and we derive the necessary estimates for the Lagrangian transformation.
In Section 3 we deal with the linearization of (\ref{system_z}). For the linear system 
we show the estimate in $W^2_p(\Omega) \times W^1_p(\Omega)$.
It is given by Theorem \ref{tw_est_main_lin}.
The first step is the energy estimate (\ref{ene}).
Then we consider the vorticity of the velocity and the Helmholtz decomposition 
to reduce the continuity equation to a sort of transport equation (\ref{trans}) that enables
us to find the bound on $\|w\|_{W^1_p(\Omega)}$. This result together with the properties of the Lam\'e system
lets us conclude Theorem \ref{tw_est_main_lin}.
In the second part of this section we apply the estimates
to solve the linear system and hence show that $T$ given by (\ref{def_T}) is well defined. 
In Section 4 we show the contraction principle for $T$. To this end
we consider the system for the difference of two solutions and write 
it in a form (\ref{system_dif}) which has a structure of (\ref{system_lin}). The contraction results from 
the estimate (\ref{est_main_lin}) and bounds on the norms on of the r.h.s. of the
system for the difference.    
At the end of Section 4 we apply the Banach fixed point theorem to solve the system (\ref{system_z})
and conclude the proof of Theorem 1. 

Let us finish this introductory part with some remarks concerning notation. By $C$ we shall denote
a constant that is controlled, but not necessarily small. $E$ shall denote a constant that can
be arbitrarily small provided the data is small enough. 
Sometimes we will write $E(\cdot)$ to underline that we need the smallness of certain quantity.
The functional spaces on $\Omega$ will be denoted without the symbol of the set, for example we
will write $W^k_p$ instead of $W^k_p(\Omega)$ for standard Sobolev spaces
of functions intergable with the $p$-th power with derivatives up to order $k$, 
$W^{1-1/p}_p(\partial \Omega)$ denotes the Slobodeckij spaces, defining regularity of traces from $W^1_p(\Omega)$, \cite{Ad}.  
Finally, we will need to consider the
density in the space $L_{\infty}(0,L;L_2(\Omega_0))$. For simplicity we denote it as $L_{\infty}(L_2)$. 
We do not use different notation for scalar and vector valued functions, while matrix valued functions are written in bolded font. 
The coordinates of a vector are denoted by $^{(\cdot)}$, i.e. $u=(u^{(1)},u^{(2)},u^{(3)})$. 

\section{Preliminaries}

In this section we introduce perturbations of the Poiseuille flow $(\bar V, \bar \rho)$
as unknown variables, what leads to the system (\ref{system}). Then we introduce a change of variables
that straightens the characteristics of the continuity equation. We obtain the system 
(\ref{system_z}). The simplified form of the continuity equation in this Lagrangian framework 
makes it possible to apply the Banach fixed point theorem to the system (\ref{system_z}).

\subsection{Reformulation of the problem}
We come back to the main system (\ref{main_system}). Since we are interested in solutions that are small 
perturbations of $(\bar V, \bar \rho)$, it is convenient to consider the perturbations as
unknown functions. For technical reasons it is better to have $u \cdot n = 0$ on the
boundary. Hence we start introducing $u_0 \in W^2_p(\Omega)$ such that   
$u_0 \cdot n|_{\Gamma} = d - \bar d$ (recall that $\bar d = \bar V \cdot n$). 
It can be found as $u_0 = \nabla \phi$ where $\phi$ solves a 
Neumann problem. We  assume that $\|d - \bar V \cdot n\|_{W^{1-1/p}_p(\Gamma)}$
is small enough for  
\begin{equation} \label{pert_pos}
V^P + u_0^{(1)}|_{\bar \Omega} \geq \theta_1 
\end{equation}
to hold for some $\theta_1>0$. In fact this is not really a restriction as
we consider small perturbations of $(\bar V, \bar \rho)$ and 
in Lemma \ref{lem_vp} we have shown that $V^P$ is separated from zero.
Now we take
\begin{equation} \label{u}
u = v - \bar V - u_0.
\end{equation}

\noindent
In particular we want our perturbed flow $v$ to have the first component also separated
from zero. This is quite natural constraint if we consider small perturbations of $\bar V$.
With the above definition of $u$ this constraint reads
\begin{equation} \label{v1_pos}
V^P + u^{(1)} + u_0^{(1)} \geq \theta_2 > 0.
\end{equation}
Next we introduce the perturbation of the density (recall that $\bar \rho \equiv 1$): 
\begin{equation} \label{w}
w = \rho - 1.
\end{equation}
%
%
%
Substituting (\ref{u}) and (\ref{w}) to (\ref{main_system})  and (\ref{d-1}) we arrive at
\begin{equation} \label{system}
\begin{array}{c}
u \cdot \nabla \bar V + V^P \, \partial_{x_1} u - \mu \Delta u - (\mu+\nu) \nabla {\rm div}\,u
+ \gamma \nabla \, w = F(u,w), \\   
V^P \, \partial_{x_1}w + {\rm div} \, u + (u+u_0) \cdot \nabla w = G(u,w), \\[6pt]
n \cdot 2 \mu {\bf D}(u) \cdot \tau_k + f (u \cdot \tau_k)|_{\Gamma} = B_k, \\
n \cdot u|_{\Gamma} = 0, \qquad 
w|_{\Gamma_{in}} = w_{in},
\end{array}
\end{equation}
where $\gamma =  \pi'(1)$ and 
$$
F(u,w) = (u + u_0) \cdot \nabla (u+u_0) - u_0 \cdot \nabla \bar V - V^P \, \partial_{x_1}u_0
- [ \pi'(w+1) -  \pi'(1)] \nabla w 
$$$$
- w \omega(f) \hat{e}_1 + (w+1) (F + \omega(f)\hat e_1)
- w (u+u_0+\bar V) \cdot \nabla (u+u_0+\bar V),
$$$$
G(u,w) = -w \, {\rm div} \, u - (w+1) \, {\rm div} \, u_0,  
$$$$
B_k = b_k - n \cdot 2 \mu {\bf D}(\bar V + u_0) \cdot \tau_k - f (\bar V + u_0) \cdot \tau_k.
$$ 

\noindent
From now on we focus on the system (\ref{system}).
Notice that $\bar V \cdot \nabla \bar V = 0$, 
hence the term 
$- w (u+u_0+\bar V) \cdot \nabla u+u_0+\bar V)$ is a higher order term and so
the form of $F$ and $G$ implies immediately the following lemma. 

\begin{lem}\label{l:2}
Let $F(u,w)$ and $G(u,w)$ be given as above, then
\begin{equation} \label{FG_lp}
\|F(u,w)\|_{L_p} + \|G(u,w)\|_{W^1_p} \leq C \big[ (\|u\|_{W^2_p} + \|w\|_{W^1_p})^2 + 
E(\|u\|_{W^2_p} + \|w\|_{W^1_p})+D_0 \big],
\end{equation} 
where $E$ denotes a small, compared to $\mu$, positive constant.
\end{lem}

\subsection{Change of variables}
With our smallness assumptions it is quite natural to solve (\ref{system})
with a fixed point argument. However, a direct application of this method fails
because of the nonlinear  term $u \cdot \nabla w$ in the hyperbolic continuity equation.
The idea to overcome this problem is to introduce a change of variables such
that this awkward term vanishes. We  look for the appropriate transformation as
$x = \psi_{u+u_0}(z)$ satisfying the identity
\begin{equation} \label{change_id}
V^P \, \partial_{z_1} = V^P \partial_{x_1} + (u+u_0) \cdot \nabla_x.
\end{equation}
In the following lemma we construct the mapping $\psi_{\bar u}$ for arbitrary
function $\bar u$ small in $W^2_p$ with vanishing normal component on the boundary $\Gamma_0$. 

\begin{lem}  \label{lem_change}
Let $\|\bar u\|_{W^2_p}$ be small enough and $\bar u \cdot n|_{\Gamma_0} = 0$. 
Then there exists a diffeomorphism $x=\psi_{\bar u}(z)$ defined on $\Omega$ 
such that $\Omega = \psi_{\bar u}(\Omega)$ 
and (\ref{change_id}) holds with $u+u_0 = \bar u$.
\end{lem}

{\it Proof.}
A key point in the proof is the fact that $V^P \geq c >0$. 
In particular we are able to divide (\ref{change_id}) by $V^P$ obtaining  
$$
\partial_{z_1} = \partial_{x_1} + \tilde u \cdot \nabla_x,
$$
where $\tilde u = \frac{\bar u}{V^P}$. Since $\bar u, V^P \in W^2_p$ we have
$\tilde u \in W^2_p$ and, since we are interested in small perturbations we 
can assume that 
\begin{equation} \label{tilde_u_small}
\|\tilde u\|_{W^2_p} << 1.
\end{equation}
Now we can follow the proof from \cite{TP2} and look for $\psi(z_1,z_2,z_3) = \psi_{z_2,z_3}(z_1)$,
where for each $(z_2,z_3) \in \Gamma_{0}$ the function $\psi_{z_2,z_3}$ is a solution to 
\begin{equation}  \label{ode}
\left\{ \begin{array}{l}
\partial_s \psi_{z_2,z_3}^{(1)} = 1, \quad
\partial_s \psi_{z_2,z_3}^{(2)} = \tilde u^2(\psi_{z_2,z_3}), \quad
\partial_s \psi_{z_2,z_3}^{(3)} = \tilde u^{(3)}(\psi_{z_2,z_3}),\\
\psi_{z_2,z_3}(0) = (0,z_2,z_3).
\end{array} \right.
\end{equation}
Due to (\ref{tilde_u_small}) we  solve (\ref{ode}) for $(z_2,z_3) \in \Gamma_{in}$
following \cite{TP2} and show that there exists a set $\Omega_{\bar u}$ such that 
$\psi(\Omega_{\bar u}) \to \Omega$ is a diffeomorphism. 
It remains to show that $\Omega_{\bar u} = \Omega$.
To this end we examine the derivatives of $\psi$. 
We have $D \psi = {\bf Id} + {\bf E}$ where
\begin{equation} \label{d_psi}
{\bf E} = \left[ \begin{array}{lcr}
\displaystyle
0 & 0 & 0 \\ 
\displaystyle         
\tilde u^2(\psi(z)) & \partial_{z_2} \int_0^{z_1} \tilde u^2(\psi(s,\bar z))ds &  \partial_{z_3} \int_0^{z_1} \tilde u^2(\psi(s,\bar z))ds \\
\displaystyle
\tilde u^{(3)}(\psi(z)) & \partial_{z_2} \int_0^{z_1} \tilde u^{(3)}(\psi(s,\bar z))ds &  \partial_{z_3} \int_0^{z_1} \tilde u^{(3)}(\psi(s,\bar z))ds    
\end{array} \right]
\end{equation}  
and $\bar z:=(z_2,z_3)$. The first row of ${\bf E}$ reduces to $0$ since 
${\bf E}_{1,i} = \partial_{z_i} \int_0^{z_1}ds = 0$ for $i=2,3$. Hence 
$$
D\psi([1,0,0]) = [1,\tilde u^2(\psi(z)),\tilde u^{(3)}(\psi(z))].  
$$
Now take $x_n \to x_0 \in \Gamma_0$ and $z_n = \phi(x_n)$. Then we have
$$
{\rm lim}_{n \to \infty} D\psi(z_n)([1,0,0]) = [1,\tilde u^2(x_0),\tilde u^{(3)}(x_0)].
$$
The latter is parallel to $\Gamma_0$ since $\tilde u \cdot n|_{\Gamma_0}=0$, hence we have 
$\phi(\Gamma_0) = \Gamma_0$ (precisely we say it in a sense of tangent spaces
since $\phi$ is defined only on $\Omega$).
To examine the behavior of tangent vectors on $\Gamma_{out}$ notice that
$$
D\psi(z)([0,\tau_1,\tau_2]) = ([0,\tilde \tau_1(\psi(z)), \tilde \tau_2(\psi(z))]),
$$
where $\tilde \tau_i$ are given by appropriate entries of $D\psi$, the important fact is that
the first coordinate vanishes. Hence if we take $x_n \to x_0$ and $z_n = \phi(x_n)$,
this time with $x_0 \in \Gamma_{out}$, then 
$$
{\rm lim}_{n \to \infty} D\psi(z_n)([0,\tau_1,\tau_2]) = [0,\tilde \tau_1(x_0),\tilde \tau_2(x_0)],
$$ 
what is parallel to $\Gamma_{out}$. It shows that $\phi(\Gamma_{out}) = \Gamma_{out}$.
Since $\phi(\Gamma_0) = \Gamma_0$ by the definition of $\phi$, we conclude that 
$\Omega_{\bar u} = \Omega$ and complete the proof.
$\square$ 

\smallskip
Now we proceed with transformation of the system (\ref{system}).
As a vector field $\bar u$ satisfying the assumptions of the Lemma we take $u+u_0$ where $u$ is the solution
of (\ref{system}). 
So far we don't know if this solution exists, our goal is to show its existence. 
Hence our approach can be regarded as working in a kind of Lagrangian coordinates \cite{MZ};
assuming that the solution exists we rewrite the system in the new variables $\psi_{u+u_0}$ induced by $u$
through (\ref{change_id}). 
Then in the new coordinates we hope to be able to apply a fixed point method to show the
existence of a solution. 
Since $(u,w)$ are perturbations that are assumed to be small, we can assume 
that $\|u\|_{W^2_p}$ is small enough that the assumptions of Lemma \ref{lem_change} are satisfied.
Hence the solution gives a well defined transformation 
$x = \psi_{u+u_0}(z)$, that we denote for simplicity by $\psi$.  
If we show in addition its uniqueness in a class of small perturbations,
then denoting $\phi = \psi^{-1}$ we have $\psi(\Omega_u) = \Omega$  
and we  come back to the original coordinates where our solution solves (\ref{system}).
Rewriting the system (\ref{system}) in coordinates $z$ yields
\begin{equation} \label{system_z}
\begin{array}{c}
u \cdot \nabla_z \bar V + (V^P \circ \psi_u) \, \partial_{z_1} u - \mu \Delta_z u - (\mu+\nu) \nabla_z {\rm div}_z\,u
+ \gamma \nabla_z \, w = \tilde F(u,w),    \\   
(V^P \circ \psi_u) \, \partial_{z_1}w + {\rm div}_z \, u = \tilde G(u,w),  \\ [6pt]
n \cdot 2 \mu {\bf D}_z(u) \cdot \tau_k + f (u \cdot \tau_k)|_{\Gamma} = B_k - n \cdot 2\mu R(u,{\bf D}) \cdot \tau_k,  \\
n \cdot u|_{\Gamma} = 0, \qquad
w|_{\Gamma_{in}} = w_{in}.
\end{array}
\end{equation}
The functions $\tilde F$ and $\tilde G$ involve commutators from the change of variables.
More precisely,
\begin{equation}  \label{tildeF}
\tilde F(u,w) = F(u,w) - u \cdot R(\bar V, \nabla) - V^P \, R(u,\partial_{x_1})
+ \mu R(u,\Delta) + (\mu + \nu) R(u, \nabla {\rm div}) - \gamma R(w,\nabla)
\end{equation}
and
\begin{equation} \label{tildeG} 
\tilde G(u,w) = G(u,w) - R(u, {\rm div}).
\end{equation}
Here the first variable in the commutator $R(\cdot,\cdot)$ denotes a function 
and the second is a differential operator. For example, 
$R(w,\nabla) := \nabla_x \, w - \nabla_z \,w$ and its $i$-th coordinate reads 
$$
R^{(i)}(w,\nabla) = \big[ \partial_{z_i}w (\phi^{(i)}_{x_i}-1) + \sum_{i \neq j} w_{z_j} \phi^{(j)}_{x_i} \big].
$$
We shall not give here precise formulas for the other commutators. 
Instead, we are now ready to give some heuristic arguments that will show
what regularity we  expect from the change of variables $\phi$.
To this end note that the commutators of the operators of order $k$ depends 
on the derivatives of $\phi$
up to order $k$. More precisely, commutators of order one contain only components of the form
$$
\nabla_z u \cdot \nabla_x \phi,  
$$
while second-order commutators contain the terms
$$
\nabla^2_z u \cdot (\nabla_x \phi)^2, \quad \nabla_z u \cdot \nabla^2_x \phi
$$
and the terms of lower order. Hence in order to find the estimates on 
$\|\tilde F(u,w)\|_{L_p}$ and $\|\tilde G(u,w)\|_{W^1_p}$ we need
$$
\nabla_x \phi \in {L_\infty}, \quad \nabla^2_x \phi \in L_p,
$$
what will be satisfied provided that $\phi \in W^2_p$
due to the imbedding $W^1_p \in L_\infty$ (recall that $p>3$).  
We should also note that, for simplicity of notation, 
the functions $F(u,w)$ and $G(u,w)$ in (\ref{tildeF}) and (\ref{tildeG})
denote exactly the same quantities as before. Hence we should keep in mind that 
they contain differential operators and so now they also contain some commutators
that we will have to control to repeat the estimate (\ref{FG_lp}).

In order to show such estimate
we should have a closer look at the derivatives of $\phi$.
With our construction of $\psi=\phi^{-1}$, it is easier to consider first the derivatives
of $\psi$ that are given by (\ref{d_psi}) with 
$\tilde u = \frac{u+u_0}{V^P}$ where $u$ is the solution to (\ref{system}),
not to be confused with $\tilde u = \frac{\bar u}{V^P}$ from the proof of Lemma \ref{lem_change}
where $\bar u$ was an arbitrary function. 
In order to find the bounds on the commutators we will need the following
smallness results for our change of variables:
\begin{lem} \label{lem_dpsi}
We have
\begin{equation} \label{est_psi1}
\sum_i\| \frac{\partial \psi^{(i)}}{\partial z_i}-1\|_{W^1_p} 
+ \sum_{i \neq j} \|\frac{\partial \psi^{(i)}}{\partial z_j}\|_{W^1_p} \leq E,
\end{equation}
\begin{equation} \label{est_psi2}
\|D^2_z \psi\|_{L_p} \leq E,
\end{equation}
\begin{equation} \label{est_phi1}
\sum_i\| \frac{\partial \phi^{(i)}}{\partial x_i}-1\|_{W^1_p} 
+ \sum_{i \neq j} \|\frac{\partial \phi^{(i)}}{\partial x_j}\|_{W^1_p} \leq E,
\end{equation}
\begin{equation} \label{est_phi2}
\|D^2_x \phi\|_{L_p} \leq E,
\end{equation}
where $E$ is  sufficiently small number comparing to $\mu$, depending on norms of perturbations
measured by $D_0$ defined in (\ref{d0}).
\end{lem}
\emph{Proof.} 
The core of the proof is in the imbedding $W^1_p \subset L_\infty$.
We start with (\ref{est_psi1}). We estimate $L_p$ norm of the entries of ${\bf E}$ (\ref{d_psi}).
These results quite directly from the form of ${\bf E}$, but needs certain
attention as ${\bf E}$ depends on $\psi$ implicitly. The entries of ${\bf E}$ without integrals
will be small provided that $\psi$ is bounded what obviously holds true. For the entries involving integrals
we  change the order of integration and derivative obtaining
$$
\partial_{z_i} \psi^{(j)} = \partial_{z_i} \int_{0}^{z_1} \tilde u(\psi(s,\bar z))ds  
= \int_{0}^{z_1} [\nabla_x \tilde u^{(j)} (\psi(s, \bar z)) \cdot \partial_{z_i} \psi(s,\bar z)]ds.
$$
By Jensen inequality we have
$$
 \int_{0}^{z_1} |\partial x_k \tilde u^{(j)} (\psi(s, \bar z))|^p 
\, |\partial_{z_i} \psi^{(k)}(s,\bar z)]|^p \,ds \leq
|z_1|^{p-1} \|\nabla_x \tilde u\|_{L_\infty}^p \int_0^{z_1} |\partial_{z_i} \psi^{(k)}(s,\bar z)]|^p \,ds.  
$$ 
Integrating the last inequality over $\Omega$ we get
\begin{equation}
\|E_{ij}\|_{L_p} \leq C (1 + \sum_{k,l} \| E_{kl} \|_{L_p}) \| \nabla_x \tilde u \|_{L_\infty}.
\end{equation}
The smallness of $\tilde u$ in $W^2_p$ and the imbedding $W^1_p \subset L_\infty$ gives (\ref{est_psi1}).  
 
To show (\ref{est_psi2}) we differentiate the entries of ${\bf E}$, let us focus on entries with integrals.
We have (we omit the sum over $k$):
$$
\partial_{z_l} \partial_{z_i}\psi^{(j)} = 
\partial_{z_l} \int_0^{z_1} [\partial_{x_k} \tilde u^{(j)}(\psi(s,\bar z)) \partial_{z_i} \psi^{(k)}(s,\bar z)] \,ds = 
$$$$ 
= \int_0^{z_1} \partial_{z_l} [\partial_{x_k} \tilde u^{(j)}(\psi(s,\bar z))] \partial_{z_i}\psi^{(k)}(s, \bar z) \,ds
+ \int_0^{z_1} \partial_{x_k} \tilde u^{(j)} \psi(s,\bar z) \partial_{z_l}\partial_{z_k} \psi^{(k)}(s,\bar z) \,ds =:I_1+I_2.
$$
Again by Jensen inequality,
\begin{equation} \label{I1}
\begin{array}{c}
|I_1|^p = \big| \int_0^{z_1} \partial_{x_m} \partial_{x_k} \tilde u^{(j)}(\psi(s,\bar z)) 
\partial_{z_l} \psi^{(m)}(s,\bar z)\partial_{z_i} \psi^{(k)}(s,\bar z) ds \big|^p \\
\leq \| \nabla_z \psi \|_{L_\infty}^2p |z_1|^{p-1} \int_0^{z_1} |\partial_{x_m} \partial_{x_k} \tilde u^{(j)} |^p  ds
\end{array}
\end{equation}
and
$$
|I_2|^p \leq \|\nabla_x \tilde u\|_{L_\infty}^p \int_0^{z_1} | \partial_{z_l}\partial_{z_k} \psi^{(k)}(s,\bar z) |^p \,ds.
$$
Integrating the above inequalities over $\Omega$ we arrive at 
\begin{equation}
\|\nabla^2 \psi\|_{L_p} \leq C(\Omega) \big[ \|\nabla_x \tilde u\|_{L_\infty} \| \nabla^2 \psi \|_{L_p}
+ \|\nabla^2_x \tilde u\|_{L_p} \big],
\end{equation}
where the term $\|\nabla_z \psi\|_{L_\infty}^2$ from (\ref{I1}) has been put into the constant.
Like in the previous estimate, the imbedding $W^1_p \subset L_\infty$
and the smallness of $\tilde u$ in $W^2_p$ yield (\ref{est_psi2}). 

To show (\ref{est_phi1}) note that the smallness of {\bf E} given by (\ref{est_psi1}) 
combined with the imbedding $W^1_p \subset L_{\infty}$
implies that ${\rm det} D \psi \geq c>0$, hence $D \psi$ is invertible and we have 
$$
D\phi = D \psi^{-1} = {\bf Id} + \tilde {\bf E},
$$  
where the elements of $\tilde {\bf E}$ can be explicitly computed in terms of ${\bf E}$.
The smallness of $\tilde u$ together with the fact that $W^1_p$ is and algebra implies
smallness in $L_p$ of the entries of $\tilde {\bf E}$, which gives (\ref{est_phi1}).
Finally (\ref{est_phi2}) is obtained by taking derivatives if $D \phi$ similarly to (\ref{est_psi2}). $\square$

\smallskip

Now we are ready to show the basic estimate on the r.h.s. of (\ref{system_z}):
\begin{lem}
Let $\tilde F$ and $\tilde G$ be defined by (\ref{tildeF}) and (\ref{tildeG}). Then we have
\begin{equation} \label{tilde_fg_lp}
\|\tilde F(u,w)\|_{L_p} + \|\tilde G(u,w)\|_{W^1_p} \leq 
C [(\|u\|_{W^2_p} + \|w\|_{W^1_p})^2 + D_0] + E \, (\|u\|_{W^2_p} + \|w\|_{W^1_p}).
\end{equation}
\end{lem}
\emph{Proof.}
The bound on $\|F(u,w)\|_{L_p} + \|G(u,w)\|_{W^1_p}$ results from (\ref{FG_lp}) in Lemma \ref{l:2}
(there are also some commutators since $F$ and $G$ involve differential operators,
but these can be estimated as follows).
We briefly justify the bounds on the commutators in $\tilde F$. To start with, the first order commutators
contain the given function $\bar V$, the functions $u$ and $w$ and the derivatives
of $\phi$, but only of the form 
\begin{equation} \label{dphi_small}
\phi^{(i)}_{x_j}, \; i \neq j \quad \textrm{and} \quad \phi^{(i)}_{x_i}-1.
\end{equation} 
Hence applying Lemma \ref{lem_dpsi} we get
\begin{equation}
\|u \cdot R(\bar V, \nabla)\|_{L_p} + \|V^P \, R(u,\partial_{x_1})\|_{L_p}
+ \|R(w, \nabla)\|_{L_p} \leq E \, (\|u\|_{W^2_p} + \|w\|_{W^1_p}).
\end{equation}
The second order commutators contain the second order derivatives of $\phi$
and the first order derivatives only of the form (\ref{dphi_small}). Hence
the application of Lemma (\ref{lem_dpsi}) yields 
\begin{equation}
\|R(u,\Delta)\|_{L_p} + \|R(u, \nabla {\rm div})\|_{L_p} 
\leq E \, \|u\|_{W^2_p}
\end{equation}
and we conclude the bound on $\|\tilde F\|_{L_p}$.
In order to estimate
$\|\tilde G(u,w)\|_{W^1_p}$ we  differentiate the commutator
$$ 
R(u, {\rm div}) = \sum u^{(i)}_{z_i} (\phi^{(i)}_{x_i} - 1) + \sum_{i \neq j} u^{(i)}_{z_j} \phi^{(j)}_{x_i}(x),   
$$
what yields
$$
\partial_{z_k} R(u, {\rm div}) = \sum_i \big[ u^{(i)}_{z_i z_k} (\phi^{(i)}_{x_i} - 1)
+ u^{(i)}_{z_i} \sum_j \phi^{(i)}_{x_i x_j} \psi^{(j)}_{z_k} \big]
$$$$
+ \sum_{i \neq j} \big[ u^{(i)}_{z_j z_k} \phi^{(j)}_{x_i} (x) 
+ u^{(i)}_{z_j} \sum_l \phi^{(j)}_{x_i x_l} \psi ^{(l)}_{z_k} \big].
$$
Applying again Lemma \ref{lem_dpsi} we get
\begin{equation} \label{rdiv}
\|R(u, {\rm div})\|_{W^1_p} \leq E \|u\|_{W^2_p}
\end{equation}
and the proof is complete. $\square$

\smallskip
\noindent
From now on we focus on the system (\ref{system_z}) instead of (\ref{system}). 
It is of crucial importance for us
that we can solve (\ref{system_z}) in the domain $\Omega$, what results from our choice of the transformation
$\psi$. 

Now we are in a position to define the operator 
\begin{equation}\label{def_T0}
T:W^2_p(\Omega) \times W^1_p(\Omega) \to W^2_p(\Omega) \times W^1_p(\Omega),
\end{equation}
to which we want to apply the Banach fixed point theorem in order to
solve the system (\ref{system_z}). Namely, we set $(u,w) = T(\bar u, \bar w)$ if
\begin{equation} \label{def_T}
\begin{array}{c}
u \cdot \nabla_z \bar V + (V^P \circ \psi_{\bar u+u_0}) \, \partial_{z_1} u - \mu \Delta_z u - (\mu+\nu) \nabla_z {\rm div}_z\,u
+ \gamma \nabla_z \, w = F(\bar u, \bar w),   \\   
((V^P \circ \psi_{\bar u+u_0})) \, \partial_{z_1}w + {\rm div}_z \, u = G(\bar u, \bar w),  \\[6pt]
n \cdot 2 \mu {\bf D}_z(u) \cdot \tau_k + f (u \cdot \tau_k)|_{\Gamma} = B_k,  \\
n \cdot u|_{\Gamma} = 0, \qquad
w|_{\Gamma_{in}} = w_{in}.
\end{array}
\end{equation} 
The point is that the term $\partial_{x_1}w + u \cdot \nabla w$ is replaced by $((V^P \circ \psi_{\bar u})+u^{(1)})\partial_{z_1}w$,
and for this term we  find a bound in $W^1_p$, what is necessary to
show the contraction property of $T$. 
From now on for simplicity we will write $\psi_{\bar u}$ instead of $\psi_{\bar u+u_0}$. 

\section{A priori bounds and solution of the linear system}
In this section we deal with the linear system:
\begin{equation} \label{system_lin}
\begin{array}{c}
u \cdot \nabla_z \bar V + (V^P \circ \psi_{\bar u}) \, \partial_{z_1} u - \mu \Delta_z u - (\mu+\nu) \nabla_z {\rm div}_z\,u
+ \gamma \nabla_z \, w = F,   \\   
((V^P \circ \psi_{\bar u})) \, \partial_{z_1}w + {\rm div}_z \, u = G,  \\[6pt]
n \cdot 2 \mu {\bf D}_z(u) \cdot \tau_k + f (u \cdot \tau_k)|_{\Gamma} = B_k,  \\
n \cdot u|_{\Gamma} = 0, \qquad 
w|_{\Gamma_{in}} = w_{in},
\end{array}
\end{equation}
with given functions $F,G,B_k,\bar u$.
Note that we take the superposition of $V^P$ and $\psi_{\bar u}$
to obtain $V^P \circ \psi_u$ for the original system and hence $V^P$ when we pass
to the original system of coordinates. The same remark concerns the function $F$, but it does
not change anything in the computations.  
We need to solve the system (\ref{system_lin}) to show that $T$ is well defined by (\ref{def_T}). 
To this end we need the appropriate estimates that we show in the first part of
this section. In the second part the linear system is solved.  

\subsection{A priori bounds}
In this section we show the estimate in $W^2_p \times W^1_p$ for the solution
of the linear system (\ref{system_lin}) in the maximal regularity regime. 
The first step is the energy estimate. It is given by the following

\begin{lem} \label{lem_ene1}
Let $(u,w)$ be a solution to the system (\ref{system_lin})
with given $(F,G,B,\bar u) \in V^* \times L_2 \times L_2(\Gamma) \times W^2_p$,
where $\bar u$ is small enough to assure 
\begin{equation} \label{vp_pos_lin}
V^P + \bar u^{(1)} \geq \theta_3>0, 
\end{equation}
for some $\theta_3>0$.
Assume that the friction $f$ is large enough on $\Gamma_{in}$ and the viscosity $\mu$
and the friction on $\Gamma_0$ satisfy (\ref{mu_large}). Then
\begin{equation} \label{ene}
\|u\|_{W^1_2} + \|w\|_{L_{\infty}(L_2)} \leq C \, [\|F\|_{V^*} + \|G\|_{L_2} + \|B\|_{L_2(\Gamma)} + \|w_{in}\|_{L_2(\Gamma_{in})}],
\end{equation}
where
\begin{equation} \label{def_V}
V = \{v \in W^1_2(\Omega): v \cdot n|_{\Gamma} = 0 \}
\end{equation}
and $V^*$ is the dual space of $V$.
\end{lem}

{\it Proof.}
We start with two basic observations.
First, since $V^P_{x_1} = 0$, by (\ref{change_id}) we have

\begin{displaymath} 
\partial_{z_1} (V^P \circ \psi_{\bar u}) = 
\frac{1}{V^P} [ \bar u^{(2)} \partial_{x_2} V^P + \bar u^{(3)} \partial_{x_3} V^P ].
\end{displaymath}
Now recall that the constant $\theta$ in (\ref{VP_positive}) is independent of the smallness of perturbation.
Hence we can assume the $\|\bar u\|_{W^2_p}$ is small compared to $\theta$ and by the imbedding $W^1_p \subset L_{\infty}$
we have

\begin{equation} \label{vp_z1}
\|\partial_{z_1} (V^P \circ \psi_{\bar u})\|_{W^1_p} \leq E(\bar u,V^P).
\end{equation}

We keep in mind that (\ref{mu_large}) and (\ref{est_VP}) hold, what implies that (\ref{vp_z1}) will be 
controlled by the viscosity, more precisely even by a constant which decreases with increasing viscosity
due to (\ref{est_VP}).
In the remaining of this section we will write $V^P$ instead of $V^P \circ \psi_{\bar u}$.
The fact that we consider the superposition does not influence the computations as we have 
(\ref{vp_z1}).
We apply the identities
\begin{equation}\label{vp_1a}
\int_{\Omega} V^P \, \partial_{z_1} |u|^2 \,dx = \frac{1}{2} \int_{\Gamma} V^P \, |u|^2 \, n^{(1)} \, d\sigma 
- \int_{\Omega} |u|^2 \, \partial_{z_1} V^P \,dx 
\end{equation}
and
\begin{equation} \label{basic_id}
\begin{array}{c}
\displaystyle
\int_{\Omega} (-\mu \Delta u - (\nu +\mu)\nabla \, {\rm div} u) \cdot v \, dx =
\int_{\Omega} 2 \mu {\bf D}(u): \nabla \,v + \nu \, {\rm div} \,u \, {\rm div} \,v \, dx - \\[4pt]
\displaystyle
\int_{\Gamma} n \cdot [2 \mu {\bf D}(u)] \cdot v \, d\sigma - \int_{\Gamma} n \cdot [\nu ({\rm div} u) {\bf{Id}} ] \cdot v \, d\sigma.
\end{array}
\end{equation}
Now we multiply (\ref{system_lin})$_1$ by $u$ and integrate.
Using the above identities, with application of the boundary conditions (\ref{system_lin})$_{3,4}$
we arrive at
\begin{equation}  \label{ene0}
\begin{array}{c}
\displaystyle
\int_{\Omega} \{ 2\mu {\bf D}(u):{\bf D}(u) + \nu |{\rm div} \,u|^2 \} \,dx
+ \int_{\Gamma_{in}} (f - \frac{V^P}{2}) |u|^2 \,d\sigma
+ \int_{\Gamma_{out}} (\frac{f}{2} + V^P) |u|^2 \,d\sigma\\[6pt]
\displaystyle
- \gamma \int_{\Omega} w \, {\rm div} \,u \,dx 
+ \int_{\Omega} (u \cdot \nabla \bar V) \cdot u \,dx = \\[6pt]
\qquad \displaystyle
= \int_{\Omega} F \cdot u \,dx + \int_{\Gamma} \{ B_1 (u \cdot \tau_1) + B_2 (u \cdot \tau_2) \} \,d\sigma
+ \int_{\Omega} |u|^2 \partial_{z_1}V^P \,dx.
\end{array}
\end{equation}
Note that by (\ref{vp_z1}) we have 
$$
| \int_{\Omega} |u|^2 \partial_{z_1}V^P \,dx | \leq E(\bar u,\nabla V^P) \|u\|_{L_2}^2.
$$
The other terms on the r.h.s. are all 'good' terms. 
The $\Gamma_{out}$ term on the l.h.s. is nonnegative
and the $\Gamma_{in}$ term will be positive for $f$ large enough.
To deal with the term $\int_{\Omega} w \, {\rm div} \,u \,dx$ we apply 
the continuity equation to express ${\rm div} \,u$ obtaining:
$$
\int_{\Omega} w {\rm div} \,u \,dx =
- \int_{\Omega} Gw \,dx - \frac{1}{2} \int_{\Omega} w^2 \, \partial_{z_1} (V^P + \bar u^{(1)}) \,dx
$$$$
+ \frac{1}{2} \int_{\Gamma_{out}} w^2 (V^P + \bar u^{(1)}) \,d\sigma
- \frac{1}{2} \int_{\Gamma_{in}} w_{in}^2 (V^P + \bar u^{(1)}) \,d\sigma.
$$
By (\ref{vp_pos_lin}) the integral over $\Gamma_{out}$ will be nonnegative, hence we have
\begin{equation} \label{est_wdivu}  
\int_{\Omega} w \, {\rm div} \,u \,dx \leq
\|G\|_{L_2}\|w\|_{L_2} + \|\partial_{z_1} (V^P + \bar u^{(1)})\|_{L_{\infty}} \|w\|_{L_2}^2 + C \|w_{in}\|^2_{L_2(\Gamma_{in})}.
\end{equation} 
To derive the $W^1_2$ - norm of $u$ we apply the Korn inequality \cite{Sol,WZ}:
\begin{equation} \label{Korn}
\int_{\Omega} [2\mu {\bf D}(u):{\bf D}(u) + \nu |{\rm div} \,u|^2] \,dx + \int_{\Gamma_{in}} f \,(u \cdot \tau)^2 \,d\sigma 
\geq C_K \| u \|_{W^1_2}^2,
\end{equation}
where $C_K=C_K(\mu,\nu,f,\Omega)$ and $C_K$ is increasing with $\mu$. A sketch of the proof 
of (\ref{Korn}) one can find in the Appendix,
note that in (\ref{Korn}) we use only information at $\Gamma_{in}$, a part of the  boundary, 
but still it is sufficient to control the whole norm of $W^1_2$.

Combining (\ref{ene0}), (\ref{est_wdivu}) and (\ref{Korn}) we get
\begin{equation} \label{ene_2}
\begin{array}{c}
\displaystyle
C_K \| u \|_{W^1_2}^2 \leq 
\|G\|_{L_2} \|w\|_{L_2} + C \|w_{in}\|^2_{L_2(\Gamma_{in})} \\
\displaystyle
+ [\|F\|_{V^*} + \|B\|_{L_2(\Gamma)}] \, \|u\|_{W^1_2} 
+ E(\bar u,\nabla V^P) [\|w\|^2_{L_2} + \|u\|_{L_2}^2]
- \int_{\Omega} (u \cdot \nabla \bar V) \cdot u \,dx.
\end{array}
\end{equation}
We have to deal with the last term on the r.h.s. It is impossible to show it has a good
sign, hence the only way is to estimate it directly with
\begin{equation} \label{ene_2a}
\big| \int_{\Omega} (u \cdot \nabla \bar V) \cdot u \,dx \big| 
\leq C_P^2 \|\nabla \bar V\|_{L_{\infty}} \|u\|_{W^1_2}^2, 
\end{equation}
where $C_P$ is the constant from the Poincar\'e inequality in $V$. Inserting the above to (\ref{ene_2})
we get
\begin{equation} \label{ene_3}
\begin{array}{c}
(C_K - C_P^2 (\|\nabla \bar V\|_{L_{\infty}} + E(\bar u,\nabla V^P)) \, \|u\|_{W^1_2}^2 \leq \\
\leq E \, \|w\|_{L_2}^2 + C ( \|G\|_{L_2}+\|w_{in}\|_{L_2(\Gamma_{in})} ) 
+ (\|F\|_{V^*}+\|B\|_{L_2(\Gamma)}) \|u\|_{W^1_2},
\end{array}
\end{equation}
where we recall that $E$ is a small constant and $C$ is a data-dependent constant, not
necessarily small. Now, $C_K$ is increasing with $\mu$, while $C_P$
does not depend on $\mu$.
Moreover, (\ref{est_VP}) implies that $E(\bar u,\nabla V^P)$ will be decreasing when $\mu$ increases. 
Finally, (\ref{nabla_pois_small}) implies that for small values of the friction $f$
on $\Gamma_0$ it is enough to assume that the viscosity is large only compared to $f|_{\Gamma_0}$
to control $E(\bar u,\nabla V^P)$.
We conclude that the constant on the l.h.s. of (\ref{ene_3}) will be positive provided that
$\mu$ satisfies (\ref{mu_large}).

Here it is a good point to emphasize the necessity of sufficient magnitude of the viscosity 
coefficient for large $f$ at $\Gamma_0$. 
This assumption is somehow natural, although
in the case $V^P=const$ it is not required \cite{TP1}. We have to control (\ref{ene_2a}), 
and largeness of dissipation may only 
give us this chance. Note that for the Dirichlet boundary condition \cite{Kw1},  
although  the constant flow is considered, such assumption is required, too.

To complete the proof of (\ref{ene}) we  find a bound on $\|w\|_{L_{\infty}(L_2)}$.
To this end we refer to the next subsection, where we solve the linear
system. Notice that $w = S(G - {\rm div} \, u)$, 
where $S$ is defined in (\ref{def_S}), and so 
\begin{equation} \label{ene_4}
\|w\|_{L_\infty(L_2)} \leq C \, (\|G\|_{L_2} + \|u\|_{W^1_2}).  
\end{equation}
Combining (\ref{ene_3}) and (\ref{ene_4}) we conclude (\ref{ene}). $\square$

\smallskip
In the next step we show higher bound on the vorticity of the velocity. To this end
we take the vorticity of (\ref{system_lin})$_1$. Denoting $\alpha = {\rm rot} \,u$ we get
\begin{equation}   \label{system_rot}
\begin{array}{lcr}
- \mu \Delta \alpha = {\rm rot}\, [F - V^P \, \partial_{z_1}u - u \cdot \nabla \bar V] & \mbox{in} & \Omega, \\
\alpha \cdot \tau_2 = (2 \chi_1 - \frac{f}{\nu}) u \cdot \tau_1 + \frac{B_1}{\nu} & \mbox{on} & \Gamma, \\
\alpha \cdot \tau_1 = (\frac{f}{\nu} - 2 \chi_2) u \cdot \tau_2 - \frac{B_2}{\nu} & \mbox{on} & \Gamma, \\
{\rm div}\,\alpha = 0 & \mbox{on} & \Gamma.
\end{array}
\end{equation}
The boundary conditions (\ref{system_rot})$_{2,3}$ are derived from differentiation of 
(\ref{system_lin})$_4$ in tangential directions and application of (\ref{system_lin})$_3$, 
see \cite{MR}, \cite{TP1}.
The above system gives the estimate (\cite{WZ}, Theorem 10.3 with $\mu=0$):
$$
\|\alpha\|_{W^1_p} \leq C \, [\|F\|_{L_p} + \|\bar V\|_{W^1_{\infty}} \|u\|_{W^1_p} + \|u\|_{W^{1-1/p}_p(\Gamma)}
+ \|B\|_{W^{1-1/p}_p(\Gamma)}]
$$$$ 
\leq C \, [\|F\|_{L_p} + \|B\|_{W^{1-1/p}+p(\Gamma)} + \|u\|_{W^1_p}].
$$
Applying the interpolation inequality (\ref{int}) to $\|u\|_{W^1_p}$ and the energy estimate we get
\begin{equation} \label{rotuw1p}
\|\alpha\|_{W^1_p} \leq C(\epsilon) \, [\|F\|_{L_p} + \|G\|_{W^1_p} + \|B\|_{W^{1-1/p}_p(\Gamma)}
+ \|w_{in}\|_{W^1_p(\Gamma_{in})}] + \epsilon \|u\|_{W^2_p}
\end{equation} 
for any $\epsilon>0$. Now consider the Helmholtz decomposition of the velocity
\begin{equation} \label{Helm}
u = \nabla \phi + A,
\end{equation}
where $\frac{\partial \phi}{\partial n}|_{\Gamma}=0$ and ${\rm div}\, A =0$.
We see that the field $A$ satisfies the system
\begin{equation}
\begin{array}{lcr}
{\rm rot}\, A = \alpha & \mbox{in} & \Omega, \\
{\rm div}\, A = 0 & \mbox{in} & \Omega, \\
A \cdot n = 0 & \mbox{on} & \Gamma.
\end{array}
\end{equation}
This is the classical rot-div system and from \cite{Sol} we have
$
\|A\|_{W^2_p} \leq C \, \|\alpha\|_{W^1_p} ,
$
what by (\ref{rotuw1p}) can be rewritten as
\begin{equation}  \label{Aw2p}
\|A\|_{W^2_p} \leq C(\epsilon) \, [ \|F\|_{L_p} + \|G\|_{W^1_p} + \|B\|_{W^{1-1/p}_p(\Gamma)}
		+ \|w_{in}\|_{W^1_p(\Gamma_{in})} ] + \epsilon \|u\|_{W^2_p}
\end{equation}
for any $\epsilon>0$.
Now we substitute the Helmholtz decomposition to (\ref{system_lin})$_1$. We get
\begin{equation}  \label{nablaH}
\begin{array}{c}
\nabla [ -(\nu + 2\mu) \Delta \phi + \gamma \, w ] = \\[7pt]
F - V^P \partial_{z_1} A + \mu \Delta A + (\nu+\mu) \nabla \,{\rm div}\,A - A \cdot \nabla \bar V 
- V^P \partial_{z_1} \nabla \phi - \nabla \phi \cdot \nabla \bar V =: \bar F.
\end{array}
\end{equation}
Since $\Delta \phi = {\rm div}\, u$, we can write
\begin{equation} \label{def_barH}
-(\nu + 2\mu) {\rm div}\,u + \gamma \, w = \bar H.
\end{equation}
We underline that we are now at the level of a priori estimates and (\ref{def_barH}) should be treated as the
definition of $\bar H$. In fact we can think of $\bar H$ as a kind of effective viscous flux like in
the theory of weak solutions to compressible Navier-Stokes equations (\cite{Fe},\cite{PL}).
Now (\ref{nablaH}) can be rewritten as
$\nabla \bar H = \bar F$. Combining the last equation with (\ref{system_lin})$_2$ we arrive at
\begin{equation} \label{trans}
\bar \gamma w + V^P \, \partial_{z_1}w = H,
\end{equation}
where $\bar \gamma = \frac{\gamma}{\nu + 2 \mu}$ and

\begin{equation} \label{H}
H = \frac{\bar H}{\nu + 2\mu}+G.
\end{equation}

The equation (\ref{trans}) makes it possible to estimate $\|w\|_{W^1_p}$ and
$\|\partial_{z_1}w\|_{W^1_p}$ in terms of $H$.
Next we can find the bound on $H$ using interpolation and the energy estimate. 
The first step is in the following lemma:
\begin{lem}
Let $w$ solve (\ref{trans}) with $H \in W^1_p$ and 
$w|_{\Gamma_{in}} = w_{in}$. Then
\begin{equation} \label{w_w1p}
\|w\|_{W^1_p} + \|\partial_{z_1}w\|_{W^1_p} \leq C \, [\|H\|_{W^1_p} + \|w_{in}\|_{W^1_p(\Gamma_{in})} ].
\end{equation}
\end{lem}
\emph{Proof.}
To estimate $\|w\|_{L_p}$ we multiply (\ref{trans}) by $|w|^{p-2}w$ and integrate.
Using the boundary conditions we get
$$
\bar \gamma \|w\|_{L_p}^p + \frac{1}{p} \int_{\Gamma_{out}} |w|^p V^P d\sigma \leq
$$$$
\leq \|H\|_{L_p} \|w\|_{L_p}^{p-1} + \frac{1}{p} \int_{\Gamma_{in}} |w|^p V^P d\sigma 
+ (\|\partial_{z_1}V^P\|_{L_{\infty}} + \|\bar u\|_{W^1_{\infty}}) \, \|w\|_{L_p}^{p-1}.
$$
The boundary term on the l.h.s. is positive and the constant in the last term on the r.h.s. is small
(note that we take only the $z_1$ derivative of $V^P$). Hence the above implies
\begin{equation} \label{w_lp}
\|w\|_{L_p} \leq C \, \big[ \|H\|_{L_p} + \|w_{in}\|_{L_p(\Gamma_{in})} \big] .
\end{equation}
In order to find a bound on $\partial_{z_i}w$
we differentiate (\ref{trans}) with respect to $z_i$. If we assume that $w \in W^1_p$ then (\ref{trans})
implies $\partial_{z_1}w \in W^1_p$, since $W^1_p$ is an algebra.
Thus we differentiate (\ref{trans}) with respect to $z_i$, 
multiply by $|\partial_{z_i}w|^{p-2} \partial_{z_i}w$ and integrate. We have

\begin{equation} \label{wzi_1}
\begin{array}{c}
\displaystyle
\int_{\Omega} (V^P+\bar u^{(1)}) |\partial_{z_i}w|^{p-2} \partial_{z_i}w \partial_{z_1}w = 
\frac{1}{p} \int_{\Omega} (V^P+\bar u^{(1)}) \partial_{z_1} |\partial_{z_i}w|^p \,dx = \\
\displaystyle
\int_{\Gamma} (V^P+\bar u^{(1)}) |\partial_{z_i}w|^p n^{(1)} d\sigma

- \int_{\Omega} |\partial_{z_i}w|^p \partial_{z_1}(V^P+\bar u^{(1)}). 
\end{array}
\end{equation}
The last term on the r.h.s. can be estimated by $E \, \|\partial_{z_i}w\|_{L_p}$ since we take only
$z_1$ derivative of $V^P$. The boundary term vanishes on $\Gamma_0$ and $\Gamma_{out}$ part will
be nonnegative. 
Hence (\ref{wzi_1}) implies 
\begin{equation} \label{wxi_2}
\bar \gamma \|\partial_{z_i}w\|_{L_p}^p \leq \|H\|_{W^1_p} \, \|\partial_{z_i}w\|_{L_p}^{p-1}
+ \frac{1}{p} \int_{\Gamma_{in}} |\partial_{z_i}w_{in}|^p (V^P + \bar u^{(1)}) \,d\sigma + E \|\partial_{z_i}w\|_{L_p}^p.
\end{equation} 
For $i=2,3$ we use the fact that $w_{in} \in W^1_p(\Gamma_{in})$ and conclude that
\begin{equation} \label{wxi_lp}
\|\partial_{z_i}w\|_{L_p} \leq C \, [\|H\|_{W^1_p} + \|w_{in}\|_{W^1_p(\Gamma_{in})}].
\end{equation}
To apply this method to $\partial_{z_1}w$ we need
some knowledge on $\partial_{z_1}w_{in}|_{\Gamma_{in}}$. To this end we can use (\ref{trans}),
which, since $V^P + \bar u^{(1)} > 0$, can be rewritten on $\Gamma_{in}$ as
$$
\partial_{z_1}w_{in} = \frac{H - \bar \gamma w_{in}}{V^P+\bar u^{(1)}}.
$$  
Hence $\|\partial_{z_1}w\|_{L_p(\Gamma_{in})} \leq C [ \|H\|_{L_p(\Gamma_{in})} + \|w_{in}\|_{W^1_p(\Gamma_{in})} ]$,
and (\ref{wxi_2}) implies (\ref{wxi_lp}) also for $i=1$.
From (\ref{w_lp}) and (\ref{wxi_lp}) we conclude
$$
\|w\|_{W^1_p} \leq C \, \big[ \|H\|_{L_p} + \|w_{in}\|_{L_p(\Gamma_{in})} \big].
$$
The bound on $\|\partial_{x_1}w\|_{W^1_p}$ results simply from the identity (\ref{trans})
and the fact that $W^1_p$ is an algebra. The proof of (\ref{w_w1p}) is complete.
$\square$

\smallskip
Now we need to find the bound on $\|H\|_{W^1_p}$, but this is straightforward.
Interpolation inequality (\ref{int}) yields
$$
\|H\|_{L_p} \leq \delta \|\nabla H\|_{L_p} + C(\delta) \, \|H\|_{L_2},
$$
for any $\delta>0$. To estimate $\|H\|_{L_2}$ we use the fact that $\bar H = -(\nu + 2\mu) {\rm div} \, u + \gamma w$
and the energy estimate (\ref{ene}). To find the bound on $\|\nabla H\|_{L_p}$ we use 
(\ref{Aw2p}), (\ref{nablaH}), then (\ref{int}) to estimate the term $\|u\|_{W^1_p}$
and finally (\ref{ene}). We obtain
\begin{equation} \label{H_w1p}
\|H\|_{W^1_p} \leq C(\delta) \, [ \|F\|_{L_p} + \|G\|_{W^1_p} + \|B\|_{W^{1-1/p}_p(\Gamma)} 
+ \|w_{in}\|_{W^1_p(\Gamma_{in})} ] + \delta \|u\|_{W^2_p}. 
\end{equation}
We are now one short step from the main result of this section.
It is given by the following
\begin{tw} \label{tw_est_main_lin}
Let $(u,w)$ be a solution to (\ref{system_lin}) with 
$(F,G,B_k,w_{in},\bar u) \in L_p \times W^1_p \times W^{1-1/p}_p(\Gamma) \times
W^1_p(\Gamma_{in}) \times W^2_p$ such that $\|\bar u\|_{W^2_p}$ is small enough
and the friction $f$ is large enough on $\Gamma_{in}$. Assume also that the viscosity
and the friction on $\Gamma_{0}$ satisfy (\ref{mu_large}). 
Then 
\begin{equation} \label{est_main_lin}
\|u\|_{W^2_p} + \|w\|_{W^1_p} + \|\partial_{z_1}w\|_{W^1_p} 
\leq C \, [ \|F\|_{L_p} + \|G\|_{W^1_p} + \|B_k\|_{W^{1-1/p}_p(\Gamma_{in})} + \|w\|_{W^1_p(\Gamma_{in})} ]. 
\end{equation}
\end{tw}
\emph{Proof.}
To close the estimate (\ref{est_main_lin}) it remains to find the bound on $\|u\|_{W^2_p}$. 
To this end notice that in particular $u$ satisfies the Lam\'e system:
\begin{eqnarray}
\begin{array}{lcr}
- \mu \Delta  u - (\nu + \mu) \nabla {\rm div}\, u =  F - \gamma \nabla  w - V^P \partial_{z_1} u - u \cdot \nabla \bar V & \mbox{in} & \Omega,\\
n\cdot 2\mu {\bf D}( u)\cdot \tau_i +f \ u \cdot \tau_i = B_i, \quad i=1,2
&\mbox{on} & \Gamma, \\
n\cdot  u = 0 & \mbox{on} & \Gamma.\\
\end{array}
\end{eqnarray}
Lemma \ref{lem_Lame} applied to the above system yields
$$
\|u\|_{W^2_p} \leq C \, [\|F\|_{L_p} + \|w\|_{W^1_p} + \|B\|_{W^{1-1/p}_p(\Gamma)} + \|u\|_{W^1_p}].
$$
Applying the interpolation inequality to the term $\|u\|_{W^1_p}$ and then the energy estimate (\ref{ene})
we get
\begin{equation} \label{est_lin_w2p_1}
\|u\|_{W^2_p} \leq C \, [\|F\|_{L_p} + \|G\|_{W^1_p} + \|w\|_{W^1_p} + \|B\|_{W^{1-1/p}_p(\Gamma)} + \|w_{in}\|_{L_2(\Gamma_{in})}].
\end{equation}
Combining this estimate with (\ref{w_w1p}) and (\ref{H_w1p}) with appropriate $\delta$ we conclude 
(\ref{est_main_lin}).
$\square$

\smallskip
Note that Theorem \ref{tw_est_main_lin} gives more than we need to solve (\ref{system_lin}), namely
the bound in $W^1_p$ of $\partial_{z_1}w$. We will use this result to show the contraction property 
for the operator $T$ in Section 4. 
\subsection{Solution of the linear system}
With the estimates that we obtained we are ready to solve the system (\ref{system_lin}).
First we define the weak solution and show its existence. Next, applying the estimate 
(\ref{est_main_lin}) we show its regularity under the appropriate regularity of the data.
\subsubsection{Weak solution}
By the weak solution to (\ref{system_lin}) we mean a couple $(u,w) \in V \times L_{\infty}(L_2)$
such that
\begin{eqnarray} \label{weak1}
\int_{\Omega} \{ v \cdot V^P \partial_{z_1} u + u \cdot \nabla \bar V) + 2 \mu {\bf D}(u) : \nabla \,v + \nu \, {\rm div}\,u \, {\rm div}\,v
- \gamma w \, {\rm div}\,v \} \,dx\nonumber\\
+ \int_{\Gamma} f (u \cdot \tau_i) \, (v \cdot \tau_i) \,d\sigma 
= \int_{\Omega} F \cdot v \,dx + \int_{\Gamma} B_i (v \cdot \tau_i) \,d\sigma
\end{eqnarray}
is satisfied $\forall \; v \in V$ and (\ref{system_lin})$_2$ is satisfied in ${\cal D'}(\Omega)$, i.e.
for all $ \; \phi \in \bar C^{\infty}(\Omega), \, \phi|_{\Gamma_{out}}=0$:
\begin{equation} \label{weak2}
-\int_{\Omega} (V^P+\bar u^{(1)}) w \partial_{z_1}\phi \,dx - \int_{\Omega} \partial_{z_1}(V^P+\bar u^{(1)}) w \phi \,dx =
\int_{\Omega} \phi (G - {\rm div}\, u) \,dx + \int_{\Gamma_{in}} (V^P + \bar u^{(1)}) w_{in} \phi \, d\sigma.
\end{equation}
To find the weak solution we apply the Galerkin method.
Hence we introduce an orthonormal basis of $\omega_k \subset V$ and finite dimensional spaces  
\mbox{$V^N = \{ \sum_{i=1}^N \alpha_i \omega_i: \; \alpha_i \in \mathbb{R} \} \subset V$}.
We look for the approximations of the velocity of the form $\displaystyle u^N = \sum_{i=1}^N c_i^N \omega_i$.  
Taking into account the continuity equation we have to define the approximations of the density
in an appropriate way. Namely, we set
$w^N = S(G^N - {\rm div} \, u^N)$, where $S: L_2(\Omega) \to L_{\infty}(L_2)$ is defined as  
\begin{equation} \label{def_S}
w=S(v) \iff \left\{ \begin{array}{lcr}
(V^P + \bar u^{(1)}) \, \partial_{z_1}w = v & \textrm{in} & {\cal D'}(\Omega), \\
w = w_{in} & \textrm{on} & \Gamma_{in}
\end{array} \right.
\end{equation}
and satisfies the estimate
\begin{equation} \label{est_S}
\|S(v)\|_{L_\infty(L_2)} \leq C \, [\|w_{in}\|_{L_2(\Gamma_{in})} + \|v\|_{L_2} ].
\end{equation} 
The construction of $S$ is quite straightforward. For a continuous $v$ we set 
\begin{equation} \label{S}
S(v)(z) = w_{in}(0,z_2,z_3)
+ \int_{0}^{z_1} \frac{v}{V^P+\bar u^{(1)}} (s,z_2,z_3) \,ds.
\end{equation}
Next we show directly the estimate (\ref{est_S}),
which enables us to extend $S$ on $L_2(\Omega)$ using a standard density argument.

Now we proceed with the Galerkin scheme. 
Taking $F = F^N$, $u = u^N = \sum_i c_i^N \, \omega_i$, $v = \omega_k, \quad k=1 \ldots N$ 
and $w = w^N = S(G^N - {\rm div}\, u^N)$ in (\ref{weak1}),
where $F^N$ and $G^N$ are orthogonal projections of $F$ and $G$ on $V^N$,
we arrive at a system of $N$ equations
\begin{equation}  \label{system_aprox}
B^N(u^N,\omega_k) = 0, \quad k = 1 \ldots N,
\end{equation}
where $B^N:V^N \to V^N$ is defined as
\begin{equation}
\begin{array}{c}
B^N(\xi^N,v^N) = \int_{\Omega} \big\{ v^N V^P \partial_{z_1}\xi^N + \xi^N \cdot \nabla \bar V + 2\mu\mathbf{D}(\xi^N) : \nabla v^N + {\rm div}\,\xi^N \, {\rm div}\,v^N \big\}\,dx \\
- \gamma \int_{\Omega} S(G^N-{\rm div}\,\xi^N) \, {\rm div}\,v^N \,dx 
+ \int_{\Gamma} [f \, (\xi^N \cdot \tau_j) - B_i] \, (v^N \cdot \tau_j) \,d\sigma - \int_{\Omega} F^N \cdot v^N \,dx.
\end{array}
\end{equation}
Now, if $u^N$ satisfies (\ref{system_aprox}) for $k = 1 \ldots N$ and $w^N = S(G^N-{\rm div}\,u^N)$, then a pair
$(u^N,w^N)$ satisfies (\ref{weak1}) - (\ref{weak2}) for $(v,\phi) \in (V^N \times \bar C^{\infty}(\Omega))$,
$\phi|_{\Gamma_{out}}=0$.
We will call such a pair an approximate solution to (\ref{weak1}) - (\ref{weak2}).
To solve the system (\ref{system_aprox}) we apply the following well-known result
(the proof can be found in \cite{Te}):
\begin{lem}  \label{lem_P}
Let $X$ be a finite dimensional Hilbert space and let $P:X \to X$ be a continuous operator
satisfying
\begin{equation}  \label{lem_P_1}
\exists M>0: \quad (P(\xi),\xi) > 0 \quad \textrm{for} \quad \|\xi\| = M.
\end{equation}
Then
$
\exists \xi^*: \quad \|\xi^*\| \leq M \quad \textrm{and} \quad P(\xi^*) = 0.
$
\end{lem} 
We define $P^N:V^N \to V^N$ as
\begin{equation}   \label{P}
P^N(\xi^N) = \sum_k B^N(\xi^N,\omega_k) \omega_k \quad \textrm{for} \quad \xi^N \in V^N.
\end{equation}
In order to apply Lemma \ref{lem_P} we show that $\big(P(\xi^N),\xi^N\big) > 0$ on some sphere in $V^N$
with radius dependent on the norms of the data.
To this end we follow the proof of the energy estimate for (\ref{system_lin}). This is in fact standard
approach in the Galerkin method: the energy estimate combined with Lemma \ref{lem_P} gives existence
of the approximate solutions, hence we skip the details here.  
Except from the existence of the approximate solution $u^N$, Lemma \ref{lem_P} gives the estimate
\begin{displaymath}  
\|u^N\|_{W^1_2} \leq C(DATA),
\end{displaymath}
which combined with (\ref{est_S}) gives
$$
\|u^N\|_{W^1_2}+\|w^N\|_{L_{\infty}(L_2)} \leq C(DATA).
$$
Thus
$$
u^N \rightharpoonup u \; {\rm in} \; W^1_2 \quad {\rm and} \quad w^N \rightharpoonup^* w \; {\rm in} \; L_{\infty}(L_2)
$$
for some $(u,w) \in W^1_2 \times L_{\infty}(L_2)$.
We easily to verify that $(u,w)$ is a weak solution.
First, passing to the limit in (\ref{weak1}) for $(u^N,w^N)$ we see that
$u$ satisfies (\ref{weak1}) with $w$. On the other hand, taking the limit in (\ref{weak2})
we verify that $w = S(G - {\rm div}\, u)$.
We conclude that $(u,w)$ satisfies (\ref{weak1}) - (\ref{weak2}), thus we have the weak solution.
To show the boundary condition on the density we can rewrite the r.h.s of (\ref{def_S}) as
\begin{equation} \label{def_S_1}
\left\{ \begin{array}{lcr}
\partial_{z_1}w = \frac{v}{V^P+\bar u^{(1)}} & \textrm{in} & {\cal D'}(\Omega), \\
w = w_{in} & \textrm{on} & \Gamma_{in},
\end{array} \right.
\end{equation}
and, treating $x_1$ as a 'time' variable, adapt Di Perna - Lions theory of transport equation (\cite{DPL})
that implies the uniqueness of solution to (\ref{def_S_1}) in the class $L_{\infty}(L_2)$
(note that this is the reason we work with weak solutions with the density in $L_\infty(L_2)$ 
instead of usual $L_2$).  
This completes the proof of existence of weak solution. 
\subsubsection{Strong solution}

The following result gives strong solution to the linear system (\ref{system_lin}) for the data
of appropriate regularity.
\begin{tw}
Let $(F,G,B_k,\bar u) \in (L_p \times W^1_p \times W^{1-1/p}_p(\Gamma) \times W^2_p)$ with 
$\|\bar u\|_{W^2_p}$ small enough. Assume further that $f$ is large enough on $\Gamma_{in}$
and $\mu$ fulfills (\ref{mu_large}). Then there exist a unique solution $(u,w) \in W^2_p \times W^1_p$
to the system (\ref{system_lin}) and the estimate (\ref{est_main_lin}) holds.
\end{tw}
{\it Proof.} To show appropriate regularity of the weak solution for the regular data
it is enough to apply the estimate (\ref{est_main_lin}) provided that we  handle the singularities
of the boundary at the junctions of the wall $\Gamma_0$ with inlet $\Gamma_{in}$ and outlet $\Gamma_{out}$. 
To this end we apply the result on the
elliptic regularity of the Lam\'e system with slip boundary conditions, Lemma \ref{lem_Lame} in the Appendix. 
Notice that we can apply this
method since we work in the fixed domain $\Omega$ due to the appropriate choice of the change
of variables. Otherwise we would end up in a free boundary problem and the solution of the linear
system would be highly nontrivial.  $\square$

\section{Contraction}

In this section we show the contraction property for the operator $T$ defined in (\ref{def_T}). 
However, first of all we  notice that Theorem 3 gives the solution
of the linear system (\ref{system_lin}) provided that $\|\bar u\|_{W^2_p}$ is small enough,
and so this constraint must hold if we want to have $T(\bar u, \bar w)$
well defined. 
We start this section with showing a stronger result, namely that 
$T:B_R \to B_R$ for some ball $B_R \in W^2_p \times W^1_p$ with
$R$ depending on the data.
This property will also be needed to show the contraction. Eventually, we prove Theorem 1, the main result of the paper.
\begin{lem} \label{lem_t1}
There exists $R>0$ depending on the size of the data measured by $D_0$ (\ref{d0}) such that, provided the data is small enough, 
$T(B_R) \subset B_R$, where $B_R$ is a ball of radius $R$
in $W^2_p \times W^1_p$.
\end{lem}
\emph{Proof.} The estimates (\ref{tilde_fg_lp}) and (\ref{est_main_lin}) imply that 
$$
\|T(\bar u,\bar w)\|_{W^2_p \times W^1_p} \leq C \, \|(\bar u,\bar w)\|_{W^2_p \times W^1_p}^2 + \delta,
$$
where $\delta < \frac{1}{4C}$, provided that the data is small enough.
For such $\delta$ we have
$$
\|(\bar u,\bar w)\|_{W^2_p \times W^1_p} < 2\delta \Rightarrow   
\|T(\bar u,\bar w)\|_{W^2_p \times W^1_p} < 2\delta. \quad
\square
$$

Now we show the contraction property for $T$. To this end
consider $(u_1,w_1) = T(\bar u_1, \bar w_1)$, $(u_2,w_2) = T(\bar u_2, \bar w_2)$.
By the definition of $T$, the difference $(u,w) := (u_1-u_2,w_1-w_2)$ satisfies the system   
\begin{equation} \label{system_dif}
\begin{array}{c}
u \cdot \nabla_z \bar V + V^P \, \partial_{z_1} u - \mu \Delta_z u - (\mu+\nu) \nabla_z {\rm div}_z\,u
+ \gamma \nabla_z \, w = 
\tilde F(\bar u_1,\bar w_1) - \tilde F(\bar u_2,\bar w_2),    \\   
(V^P+\bar u_2^{(1)}) \, \partial_{z_1}w + {\rm div}_z \, u = 
\tilde G(\bar u_1,\bar w_1) - \tilde G(\bar u_2,\bar w_2) + (\bar u_1 - \bar u_2) \partial_{z_1}w_1, \\[6pt]
n \cdot 2 \mu {\bf D}_z(u) \cdot \tau_k + f (u \cdot \tau_k)|_{\Gamma} = n \cdot 2\mu [R(u_2,{\bf D})-R(u_1,{\bf D})] \cdot \tau_k,  \\
n \cdot u|_{\Gamma} = 0, \qquad
w|_{\Gamma_{in}} = 0.
\end{array}
\end{equation}
Hence in order to show the contraction principle for $T$ we can apply (\ref{est_main_lin})
provided that we have good bounds on the r.h.s. of (\ref{system_dif}). 
A result we need is given by the following 
\begin{lem}
We have
\begin{equation} \label{est_dif}
\begin{array}{c}
\|\tilde F(\bar u_1,\bar w_1) - \tilde F(\bar u_2,\bar w_2)\|_{L_p} +
\|\tilde G(\bar u_1,\bar w_1) - \tilde G(\bar u_2,\bar w_2)\|_{W^1_p} \\
+ \|(\bar u_1 - \bar u_2) \partial_{z_1}w_1\|_{W^1_p} +  
\|n \cdot 2\mu [R(u_2,{\bf D})-R(u_1,{\bf D})] \cdot \tau_k\|_{W^{1-1/p}_p(\Gamma)}  \\
\leq E \, [\|u_1-u_2\|_{W^2_p} + \|w_1-w_2\|_{W^1_p}].
\end{array}
\end{equation}
\end{lem}
\emph{Proof.} 
Denote $(\bar u_i, \bar w_i)$ by $(u_i, w_i)$ for simplicity. 
We have
$$
\tilde F(u_1,w_1) - \tilde F(u_2,w_2) = u_1 \cdot \nabla_x u_1 - u_2 \cdot \nabla_x u_2 
+ (u_1 - u_2) \cdot \nabla_x u_0 + u_0 \cdot \nabla_x (u_1-u_2)
$$$$
+ \delta \pi'(w_1) \nabla_x w_1 - \delta \pi'(w_2) \nabla_x w_2 + (w_1-w_2)F  +\tilde R,
$$
where  
$$
\delta \pi'(w_i) := \pi'(w_i+1) - \pi'(w_i)
$$
and $\tilde R$ denotes all the differences between corresponding commutators in $\tilde F$.
We can afford using such abbreviation and estimate $\tilde R$ without any additional computation 
if we just notice that the commutators are linear with respect to the functions and 
hence if certain estimate in terms of the function holds for a commutator, then the same estimate
in terms of the difference hold for the difference of the commutators, for example
$$
|R(u,\partial_{x_1})| \leq E \|u\|_{W^1_p} \Rightarrow
|R(u_1,\partial_{x_1}) - R(u_2,\partial_{x_1})| \leq E \, \|u_1 - u_2\|_{W^1_p}.
$$
Applying this reasoning to all the commutators in $\tilde F$ we conclude that 
\begin{equation} \label{difF_1}
\|\tilde R\|_{L_p} \leq E(\phi) \big( \|u_1 - u_2\|_{W^2_p} + \|w_1 - w_2\|_{W^1_p} \big) .
\end{equation}
We estimate the remaining parts. Obviously we have
\begin{equation} \label{difF_2}
\|(u_1 - u_2) \cdot \nabla_x u_0\|_{L_p} + \|u_0 \cdot \nabla_x (u_1-u_2)\|_{L_p}
\leq E(u_0,\phi) \|u_1 - u_2\|_{W^2_p},
\end{equation}   
where $E(\cdot,\cdot)$ depends also on $\phi$ since we have commutators as the gradients are w.r.t. $x$.
A little bit closer examination shows that the $i$-th coordinate
$$
\big( u_1 \cdot \nabla_x u_1 - u_2 \cdot \nabla_x u_2 \big)^{(i)} = 
\big( u_1 \cdot \nabla_x (u_1 - u_2) + (u_1-u_2) \cdot \nabla_x u_2 \big)^{(i)} =
$$$$ 
= \sum_j \big[ (u_1 - u_2)^{(j)} \partial_{z_j}u^{(i)}_2 + R(u_2^{(i)},\partial_{x_j}) \big]
+\sum_j \big[ u_1^{(j)} \partial_{z_j}(u_1-u_2)^{(i)} + R((u_1-u_2)^{(i)},\partial_{x_j}) \big],
$$
and so by a direct computation we get
\begin{equation} \label{difF_3}
\|u_1 \cdot \nabla_x u_1 - u_2 \cdot \nabla_x u_2\|_{L_p} 
\leq E(\|u_i\|_{W^2_p},\phi) \|u_1-u_2\|_{W^2_p}.
\end{equation}
It remains to estimate 
$$
\delta \pi'(w_1) \nabla_x w_1 - \delta \pi'(w_2) \nabla_x w_2 = 
\delta \pi'(w_1) \nabla_x (w_1-w_2) + [\delta \pi'(w_1) - \delta \pi'(w_2)] \nabla_x w_2 = 
$$$$  
= \delta \pi'(w_1) [\nabla_z (w_1-w_2) + R(w_1-w_2,\nabla)]
+ [\delta \pi'(w_1) - \delta \pi'(w_2)] [\nabla_z w_2 + R(w_2,\nabla)].
$$
It follows easily that
\begin{equation} \label{difF_4}
\|\delta \pi'(w_1) \nabla_x w_1 - \delta \pi'(w_2) \nabla_x w_2\|_{L_p}
\leq E(\|w_1\|_{W^1_p},\|w_2\|_{W^1_p},\phi) \|w_1 - w_2\|_{W^1_p}.
\end{equation}
Combining (\ref{difF_1}), (\ref{difF_2}) , (\ref{difF_3}) and (\ref{difF_4})
we conclude 
\begin{equation} \label{est_difF}
\|\tilde F(u_1,w_1) - \tilde F(u_2,w_2)\|_{L_p} \leq E(\|w_i\|_{W^1_p},\|u_i\|_{W^2_p},\phi) 
\big[ \|u_1-u_2\|_{W^2_p} + \|w_1-w_2\|_{W^1_p} \big] .
\end{equation}
Now we estimate the difference in $G$. We have
$$
\tilde G(u_1,w_1) - \tilde G(u_2,w_2) =
$$$$
=(w_2-w_1) {\rm div}u_0 + w_1 \, {\rm div}(u_1-u_2) + {\rm div} u_2 (w_1-w_2) 
+ R(u_2,{\rm div}) - R(u_1,{\rm div}).
$$
The first term
$$
\|(w_2-w_1) {\rm div}u_0\|_{W^1_p} \leq 
\|\nabla(w_2-w_1) {\rm div}u_0\|_{L_p} + \|(w_2-w_1) \nabla^2 u_0\|_{L_p} \leq 
$$$$
\leq \|{\rm div}u_0\|_{L_\infty} \|w_1-w_2\|_{W^1_p} + \|u_0\|_{W^2_p} \|w_1-w_2\|_{L_\infty} \leq
C \|u_0\|_{W^2_p} \|w_1-w_2\|_{W^1_p}. 
$$
The second
$$
\|w_1 \, {\rm div}(u_1-u_2)\|_{W^1_p} \leq
\|(\nabla w_1) \, {\rm div}(u_1-u_2)\|_{L_p} + \|w_1 \, \nabla^2 (u_1-u_2)\|_{L_p} \leq
C \, \|w_1\|_{W^1_p} \|u_1 - u_2\|_{W^2_p} .
$$
Similarly we show
$$
\|{\rm div} u_2 (w_1-w_2)\|_{W^1_p} \leq C \|u_2\|_{W^2_p} \, \|w_1-w_2\|_{W^1_p}.
$$
Now we have to estimate the difference of the commutators in $W^1_p$. 
It turns out to be straightforward as we have 
$$
R(u_1,{\rm div}) - R(u_2,{\rm div}) = \sum_i \partial_{z_i}(u_1 - u_2)^{(i)} \phi^{(i)}_{x_i}
+ \sum_{i \neq j} \partial_{z_j}(u^{(1)}-u^2)^{(1)} \phi^{(j)}_{x_i}
$$ 
and so, identically as in (\ref{rdiv}), we show that 
$$
\|R(u_1,{\rm div}) - R(u_2,{\rm div})\|_{W^1_p} \leq E(\phi) \|(u_1-u_2)\|_{W^2_p}.
$$
Combining the above results we conclude 
\begin{equation} \label{est_difG}
\|\tilde G(u_1,w_1) - \tilde G(u_2,w_2)\|_{W^1_p} \leq
\big[ E(\phi) + C \, (\|w_1\|_{W^1_p} + \|u_2\|_{W^2_p}) \big] (\|u_1-u_2\|_{W^2_p} + \|w_1 - w_2\|_{W^1_p}).
\end{equation}
To treat the last term of the r.h.s. of $(\ref{system_dif})_2$ we observe that
\begin{equation}\label{3term}
 \| (\bar u_1 -\bar u_2) \partial_{z_1}w_1\|_{W^1_p} \leq C\|\partial_{z_1}w_1\|_{W^1_p}\|\bar u_1 -\bar u_2\|_{W^2_p}.
\end{equation}
It remains to estimate the boundary terms 
$n \cdot [R(u_1,{\bf D})-R(u_2,{\bf D})] \cdot \tau_k$.
To this end it is enough to notice that
$$ 
\{R(u_1,{\bf D})-R(u_2,{\bf D})\}_{i,j} =
R((u_1-u_2)^{(i)}, \partial_{x_j}) - R((u_1-u_2)^{(j)} , \partial_{x_i}) =: R_{i,j}.
$$
Applying the trace theorem and repeating the proof of (\ref{est_difG})
we can show that
$$
\|R_{i,j}\|_{W^{1-1/p}_p(\Gamma)} \leq C \, \|R_{i,j}\|_{W^1_p} \leq E(\phi) \, \|u_1-u_2\|_{W^2_p},
$$
and so 
\begin{equation} \label{est_difRD}
\|n \cdot [R(u_1,{\bf D})-R(u_2,{\bf D})] \cdot \tau_k\|_{W^{1-1/p}_p(\Gamma)} 
\leq E(\phi) \, \|u_1-u_2\|_{W^2_p}. 
\end{equation}
Combining (\ref{est_difF}), (\ref{est_difG}), (\ref{3term}) and (\ref{est_difRD})
we conclude (\ref{est_dif}). $\square$

\smallskip
{\it Proof of Theorem 1.}
With the results of the previous section we can apply the Banach fixed point theorem
to the operator $T$. It gives existence of a unique fixed point within the ball 
$B(0,R) \in W^2_p \times W^1_p$. By Lemma \ref{lem_t1} we have 
$R = R(D_0)$ where $D_0$ is defined in (\ref{d0}). By the definition of $T$, the
fixed point $(u,w)$ solves the system (\ref{system_z}). 

Now we recall Section 2 and conclude that the original coordinate system is $x = \psi_u(z)$,
and in the $x$ variable our solution satisfies the system (\ref{system}).
It follows that $v=\bar v + u+ u_0$ and $\rho = w+1$ solves (\ref{main_system}) and the estimate
(\ref{est_main}) holds. Theorem 1 is proved.

\section{Appendix}

{\it Proof of Lemma \ref{lem_vp}.} As explained in Section 1, we assume the pressure 
of the form $\bar \Pi = \omega(f) \, x_1$. Then
on each $x_1$ - cut of $\Omega$ (i.e. on each set $\Omega_0 \times \{x_1\}$)
$V^P$ can be found as a solution to the elliptic problem (\ref{eqn_pois}) which we recall here:
\begin{equation} \label{eqn_pois1}
\begin{array}{c}
\mu \Delta v = \bar \Pi_{x_1} = \omega(f) < 0 \mbox{ ~~ in~} \Omega_0, \\
\mu \frac{\partial v}{\partial n} + f \, v = 0 \mbox{ ~~ on ~} \partial \Omega_0.
\end{array}
\end{equation}
Testing (\ref{eqn_pois}) with $v_{-} = v \, \chi_{v<0}$ 
we get
\begin{equation} \label{lem_vp_1}
- \int_{\Omega} \mu |\nabla v_{-}|^2 - \int_{\Omega} f \, v_{-}^2 = \int_{\Omega} \omega(f) \, v_{-} \geq 0.
\end{equation}
The last inequality results from nonpositivity 
of $v_{-}$ and the direction of the flow which implies $\omega(f)<0$ (recall the remark after
(\ref{eqn_pois})) 
Clearly (\ref{lem_vp_1}) implies $v \geq 0$. We want to show sharp inequality. 
To this end consider $\bar v$ satisfying
\begin{displaymath}
\begin{array}{c}
\Delta \bar v = \omega(f)/\mu  < 0\mbox{ ~~ in~} \Omega , \qquad 
\bar v|_{\Gamma} = 0 \mbox{ ~~ at ~} \partial \Omega.
\end{array}
\end{displaymath}
We have $\bar v \geq 0$, and, by the maximum principle applied to $v - \bar v$ we get
${\rm inf}_{\Omega} \, v = {\rm inf}_{\Gamma} \, v$.     
Assume that ${\rm inf}_{\Gamma}v = v(x_0) = 0$ for some $x_0 \in \Gamma$.
Then, since $f \geq 0$ and $\frac{\partial \bar v}{\partial n} \leq 0$, we must have 
$$
\frac{\partial(v - \bar v)}{\partial n} (x_0) = 
- \frac{\partial \bar v}{\partial n}(x_0) - f \, v(x_0) \geq 0.
$$  
But since $v(x_0) = {\rm inf}_{\Omega} (v - \bar v)$, by the Hopf Lemma we must have 
$\frac{\partial(v - \bar v)}{\partial n} (x_0) < 0$. 
The application of of Hopf lemma is possible since  
$\Omega_0$ is a $C^2$ subset of $\mathbb{R}^2$, hence $v$ is a classical solution to (\ref{eqn_pois1}),
 i.e. $v\in C^2(\Omega)\cap C^1(\overline{\Omega})$.
We conclude that $v \geq \theta > 0$ on $\Gamma$, hence $v \geq \theta$ in $\bar \Omega$.
In particular, the linearity of (\ref{eqn_pois1}) gives (\ref{est_VP}) 
where $\bar \omega$ is continuous since $\omega$ is continuous (recall the discussion after the formulation
of Lemma \ref{lem_vp}).
This completes the proof.
$\square$ 

\begin{lem} \label{lem_Lame} (Lam\'e system with slip boundary conditions).
Let $\mu>0$, $\nu + 2\mu>0$, $F \in L_p(\Omega)$ and $B \in W^{1-1/p}_p(\Gamma)$. 
Then there exists $u \in W^2_p(\Omega)$ solving 
\begin{equation} \label{Lam\'e}
\begin{array}{lcr}
-\mu \Delta  u - (\nu + \mu) \nabla {\rm div}\, u = F & \mbox{in} & \Omega, \\
n\cdot 2\mu {\bf D}( u)\cdot \tau +f \ u \cdot \tau = B, \quad i=1,2
&\mbox{on} & \Gamma, \\
n\cdot  u = 0 & \mbox{on} & \Gamma.\\
\end{array}
\end{equation}
Moreover, the following estimate holds: 
\begin{equation} \label{est_lame_3d}
\|u\|_{W^2_p} \leq C [\|F\|_{L_p} + \|B\|_{W^{1-1/p}_p(\Gamma)}].
\end{equation}
\end{lem}
\emph{Proof.} 
Under the assumptions on $\mu$ and $\nu$ (\ref{Lam\'e}) is elliptic so we easily
get a weak solution. 
The only problem we encounter showing regularity of the weak solution 
under appropriate regularity of the data are the singularities of the boundary
on the junctions of $\Gamma_0$ with $\Gamma_{in}$ and $\Gamma_{out}$.
These can be dealt with using symmetry arguments. But first we have to reformulate slightly the system (\ref{Lam\'e}).
Having the weak solution, and the bound in $W^1_2$ for the velocity, we are able to consider only the case $f\equiv 0$ and
$B\equiv 0$. The general case will be a consequence of this particular one. Assume that $\Gamma_{in} \subset \{x_1=0\}$,
then we define the operator $E^v_{as}$ extending a vector field defined
for $\{x: x_1 \geq 0\}$ on the whole space as 
\begin{equation}
E^v_{as}(u)(x) = \left\{ \begin{array}{c}
u(x), \quad x_1 \geq 0, \\
\tilde u(\tilde x), \quad x_1 < 0,
\end{array} \right.
\end{equation}        
where $\tilde x = (-x_1,x_2,x_3)$ and $\tilde u(\tilde x) = [-u^{(1)}(x),u^2(x),u^{(3)}(x)]$.
Then we have 
\begin{equation}
\Delta E^v_{as}(v) + \nabla {\rm div} E^v_{as}(v) = \Delta v + \nabla {\rm div} v,
\end{equation}
and on the plane ${x_1=0}$ the extension $E^v_{as}$ preserves the slip boundary conditions for $f\equiv 0$ and $B\equiv 0$, since
$$
n\cdot {\bf D}( E^v_{as}(v) ) \cdot \tau_i|_{x_1=0} = \partial_{\tau_i}E^v_{as}(v)^{(1)}_+ \partial_{x_1}E^v_{as}(v)^{(i)}=0
$$
and $n \cdot E^v_{as}(v)|_{x_1=0} = 0$, $\partial_{x_1} E^v_{as}(v)^{(i)}=0$ for $i=2,3$.
Now we are allowed to localize the equations in a vicinity of $\Gamma_{in}$ obtaining a system in a smooth domain. 
Then the theory  \cite{ADN1,ADN2} gives   the full regularity of $E^v_{as}(v)$
in neighborhood of the junctions.
Application of analogous antisymmetric extension on $\Gamma_{out}$
completes the proof. $\square$

\begin{lem} (interpolation inequality): \label{lem_int} \\
For $f \in W^1_p(\Omega), \quad p>3$:
\begin{equation}  \label{int}
\|f\|_{L_p} \leq \epsilon \|\nabla f\|_{L_p} + C(\epsilon,p,\Omega) \|f\|_{L_2}.
\end{equation}
\end{lem}
\emph{Proof.} The interpolation inequality in $L_p$ (\cite{Ad}, Theorem 2.11)
and the imbedding $W^1_p \subset L_{\infty}$ for $p>3$ yields 
\begin{displaymath}  
\|f\|_{L_p} \leq C(p) \, \|f\|_{L_\infty}^{\theta} \, \|f\|_{L_2}^{1-\theta} \leq
C(p) (\|f\|_{L_p} + \|\nabla f\|_{L_p})^{\theta} \, \|f\|_{L_2}^{1-\theta},
\end{displaymath}
what entails (\ref{int}) after application of the Cauchy inequality. $\square$
\smallskip

\begin{lem}\label{l:korn}
 Let $\Omega$ be bounded with sufficiently smooth boundary and $\Gamma_{part}$ be an open regular  subset of $\partial \Omega$. Then 
\begin{equation}\label{k0}
 \|u\|_{W^1_2} \leq C(\|\mathbf D(u)\|_{L_2} +\|u|_{\Gamma_{part}}\|_{L_2(\Gamma_{part})})
\end{equation}
for $u \in W^1_2(\Omega)$.
\end{lem}

{\it Proof.} The known result based on properties of the kernel of  $\mathbf D(\cdot)$ \cite{LadSol,Zaj} yields
\begin{equation}\label{k1}
 \|u\|_{W^1_2} \leq C_1\|\mathbf D(u)\|_{L_2} + C_2\|u\|_{L_2}.
\end{equation}
In order to obtain (\ref{k0}) we shall prove that 
\begin{equation}\label{k2}
 \|u\|_{L_2} \leq \frac{1}{2C_2} \|\nabla  u \|_{L_2} + M(\|\mathbf D(u)\|_{L_2} +\|u|_{\Gamma_{part}}\|_{L_2(\Gamma_{part})}).
\end{equation}
Compactness argument and features of $\mathbf D$ implies (\ref{k2}). We use that fact that the only solution to the system $\mathbf D(u^*)=0,
u|_{\Gamma_{part}}=0$ is $u^*\equiv 0$. To show it we use the fact that if $\mathbf D(u)=0$ then 

\begin{equation}
 u=A \left( \begin{array}{c}
x_2 \\
-x_1 \\
0
\end{array}
\right)
+ B\left( \begin{array}{c}
x_3 \\
0 \\
-x_1
\end{array}
\right)
+C \left( \begin{array}{c}
0 \\
x_3 \\
-x_2
\end{array}
\right)
+ \left( \begin{array}{c}
d_1 \\
d_2 \\
d_3
\end{array}
\right).
\end{equation}
Hence $u$ is an affine map and it is easy to verify 
verify that if at least one of the coefficients $A,B,C$ is nonzero
then the rank of matrix of $u$ is $2$, hence Ker $u$ is a line. On the other hand,  
$\Gamma_{part}$ is a two dimensional submanifold as it is an open, regular subset of $\partial \Omega$. 
But $\Gamma_{part} \subset {\rm Ker}u$, hence we conclude that $A=B=C=0$, hence $d_i=0$.  \qed 

\medskip 

{\bf Acknowledgements.} This work has been supported  partly by the Ministry of Science grant N N201 547438 
and by Foundation for Polish Science in fr. EU European Regional Development Funds (OPIE 2007-2013).


\begin{thebibliography}{99}

\footnotesize
\bibitem{Ad} Adams, R.;  Fournier, J. \emph{Sobolev spaces}\/, 2nd ed., Elsevier, Amsterdam, 2003
\bibitem{ADN1} Agmon, S.; Douglis, A.; Nirenberg, L. \emph{Estimates near the boundary for solutions of elliptic
			partial differential equations satisfying general boundary conditions I}\/,
			Comm.Pure Appl.Math. 12 (1959), 623-727

\bibitem{ADN2} Agmon, S.; Douglis, A.; Nirenberg, L. \emph{Estimates near the boundary for solutions of elliptic
			partial differential equations satisfying general boundary conditions II}\/,
			Comm.Pure Appl.Math. 17 (1964), 35-92

\bibitem{BdV} Beirao da Veiga, H. \emph{An $L^p$-Theory for the n-Dimensional, Stationary, Compressible Navier-Stokes Equations, 
			and the Incompressible Limit for Compressible Fluids. The Equilibrium Solutions}, Comm.Math.Phys. 109 (1987), 229-248			
		
\bibitem{PC} Constantin, P.
\emph{An Eulerian-Lagrangian approach for incompressible fluids: local theory.}
J. Amer. Math. Soc. 14,2 (2001), 263-278

\bibitem{DPL} R.J.DiPerna, P.L.Lions, \emph{Ordinary differential equations, transport theory and Sobolev spaces}\/,
			Invent.math. 98 (1989), 511-547
	
\bibitem{Fe} Feireisl, E. \emph{Dynamics of viscous compressible fluids}\/,
			Oxford Lecture Series in Mathematics and its Applications, 26.
			Oxford University Press, Oxford, 2004
			
\bibitem{Ga} Galdi, G.P. \emph{An Introduction to the mathematical theory of the Navier-Stokes Equations}\/,
				Vol.I, Springer-Verlag, New York, 1994
			
\bibitem{Kw1} Kellogg, R. B.; Kweon, J. R. \emph{Compressible Navier-Stokes equations in a bounded domain
				with inflow boundary condition}\/, SIAM J.Math.Anal. 28,1(1997), 94-108
				
\bibitem{Kw2} Kellogg, R. B.; Kweon, J. R. \emph{Smooth Solution of the Compressible Navier-Stokes Equations
				in an Unbounded Domain with Inflow Boundary Condition}\/,
				J.Math.Anal. and App. 220 (1998), 657-675


\bibitem{LadSol} Ladyzhenskaja, O. A.; Solonnikov, V. A. Some problems of vector analysis, and generalized formulations of boundary value problems for the Navier-Stokes equation. (Russian) Boundary value problems of mathematical physics and related questions in the theory of functions, 9.  Zap. Nau~n. Sem. Leningrad. Otdel. Mat. Inst. Steklov. (LOMI)  59  (1976), 81--116, 256.

\bibitem{PL} Lions, P.L. \emph{Mathematical topics in fluid mechanics. Vol. 2. Compressible models}\/,
			Oxford Lecture Series in Mathematics and its Applications, 10.
			Oxford Science Publications. The Clarendon Press, Oxford University Press, New York, 1998		
			
\bibitem{Mu} Mucha, P. B. \emph{The Cauchy problem for the compressible Navier-Stokes equations in the $L_p$-framework.} Nonlinear Anal. 52,4 (2003), 1379-1392.
	
\bibitem{PM1} Mucha, P. B. \emph{The Navier-Stokes equations and the maximum principle},
			Int. Math. Res. Not. 67 (2004), 3585-3605 

\bibitem{MR} Mucha, P .B.; Rautmann, R. \emph{Convergence of Rothe's scheme for the Navier-Stokes equations with slip conditions in 2D domains}\/,
			ZAMM Z. Angew. Math. Mech. 86,9 (2006), 691-701
			
\bibitem{MZ} Mucha, P. B.; Zaj\c aczkowski, W. On local existence of solutions of the free boundary problem for an incompressible viscous self-gravitating fluid motion.  Appl. Math. (Warsaw)  27,3 (2000), 319-333. 

\bibitem{NoPa} Novotn\'y, A.; Padula, M. \emph{$L^p$ approach to steady flows of viscous compressible fluids in exterior
			domains} Arch.Rat.Mech.Anal. 126(1994), 243-297  	

\bibitem{NoS} Novotn\'y, A.; Straskraba, I. \emph{An Introduction to the Mathematical Theory of Compressible Flows}\/,
			Oxford Science Publications, Oxford 2004			
\bibitem{TP1} Piasecki, T. \emph{Steady compressible Navier-Stokes flow in a square},
			J.Math.Anal.Appl. 357 (2009), 447-467
\bibitem{TP2} Piasecki, T. \emph{On an inhomogeneous slip-inflow boundary value problem for a steady
			flow of a viscous compressible fluid in a cylindrical domain}, 
			Journal of Differential Equations 248 (2010), 2171-2198

\bibitem{PRS1} Plotnikov, P. I.; Ruban, E. V.; Sokolowski, J. \emph{Inhomogeneous boundary value problems for
			compressible Navier-Stokes Equations: well-posedness and sensitivity analysis}\/,
			SIAM J.Math.Anal. 40,3 (2008), 1152-1200
			
\bibitem{PRS2} Plotnikov, P.I.; Ruban, E.V.; Sokolowski, J.	\emph{Inhomogeneous boundary value problems
			for compressible Navier-Stokes and transport equations}\/,
			J.Math.Pures Appl. 92,2 (2009), 113-162
			
\bibitem{Sol} Solonnikov, V. A. \emph{Overdetermined elliptic boundary value problems},
			Zap.Nauch.Sem.LOMI 21 (1971), 112-158	

\bibitem{Sol2} Solonnikov, V. A. \emph{On the nonstationary motion of an isolated volume of a viscous incompressible fluid}, Izv. Akad. Nauk SSSR 51 (1987), 1065-1087 (in Russian). 
\bibitem{Te} R.Temam, \emph{Navier Stokes Equations}\/,	North-Holland, Amsterdam, 1977.
\bibitem{VZ} Valli, A.; Zaj\c aczkowski, W. M. \emph{Navier-Stokes equations for compressible fluids: global existence and qualitative properties of the solutions in the general case}. 
Comm. Math. Phys. 103,2 (1986), 259-296.

\bibitem{WZ} Zaj\c aczkowski, W. M. \emph{Existence and regularity of solutions of some elliptic systems in domains
with edges},\/ Dissertationes Math., 274(1989), 95 pp.

\bibitem{Zaj} Zaj\c aczkowski, W.M. \emph{On nonstationary motion of a compressible viscous fluid bounded by a free surface}, Dissertationes Math. 324 (1993).


\end{thebibliography}
\end{document}